\documentclass[12pt]{article}
\usepackage[cp1251]{inputenc}
\usepackage[T2A]{fontenc}
\usepackage[english]{babel}

\def\eps{\varepsilon}

\newcounter{num}[section]

\newcommand{\Th}{\refstepcounter{num}
{\bf Theorem \arabic{section}.\arabic{num} }}

\newcommand{\Lemma}{\refstepcounter{num}
{\bf Lemma \arabic{section}.\arabic{num} }}

\newcommand{\Pred}{\refstepcounter{num}
{\bf Proposition \arabic{section}.\arabic{num} }}

\newcommand{\Cor}{\refstepcounter{num}
{\bf Corollary \arabic{section}.\arabic{num} }}

\newcommand{\Note}{\refstepcounter{num}
{\it Note \arabic{section}.\arabic{num} }}

\newcommand{\St}{\refstepcounter{num}
{\bf Statement \arabic{section}.\arabic{num} }}

\newcommand{\Def}{\refstepcounter{num}
{\it Definition \arabic{section}.\arabic{num} }}

\newcommand{\Proof}{{\bf Proof. }}

\def\v{\vec}

\def\a{\alpha}
\def\d{\delta}
\def\m{\times}

\def\F{\widehat}

\def\L{\Lambda}
\def\l{\lambda}

\def\k{\kappa}

\def\t{\tilde}

\def\o{\overline}

\def\ind{{\rm ind}}

\author{I.D. Shkredov }
\title{On a Generalization of Szemeredi's Theorem.}
\date{}
\begin{document}
\maketitle

\refstepcounter{section}

\begin{center}
\textsc{ \arabic{section}. Introduction.}
\end{center}

 Let $N$ be a natural number. We set
$$
  a_k(N) = \frac{1}{N} \max \{ |A| ~:~ A \subseteq [1,N],
$$
$$
  A \mbox{ contains no arithmetic progressions of length  } k
   \},
$$
where $|A|$ denotes the cardinality of $A$.
In \cite{EaT}, Erdos and Turan conjectured that any set of positive density
contains an arithmetic progression of given length.
In other words, they supposed that, for any $k\ge 3$
\begin{equation}
  a_k(N) \to 0 \mbox{ as } N\to \infty
\label{Sz_lab}
\end{equation}
Clearly, this conjecture implies van der Waerden's theorem \cite{Wdv}.

   In the simplest case of $k=3$ conjecture (\ref{Sz_lab}) was proven
 in \cite{Rt} by K.F. Roth, who applied the Hardy -- Littlewood method to show that
$$
  a_3(N) \ll \frac{1}{\log \log N}.
$$
   At present, the best upper bound for $a_3(N)$
   is due to J. Bourgain.
He proved that
\begin{equation}
 a_3(N) \ll \sqrt{ \frac{\log \log N}{ \log N} }.
\end{equation}
   For an arbitrary $k$ conjecture (\ref{Sz_lab}) was proven by
E. Szemeredi \cite{Sz} in 1975. Szemeredi's proof uses difficult combinatorial arguments.

An alternative proof was suggested by Furstenberg in \cite{Fu}.
His approach uses the methods of ergodic theory.
Furstenberg showed that Szemeredi's theorem is equivalent to the
multiple recurrence of almost all points in any dynamical system.
Here we formulate his theorem in the case of metric spaces:
\\
\Th {\it Let $X$ be a metric space with metric
$d(\cdot,\cdot)$ and Borel sigma--algebra of measurable sets
$\Phi$. Let $T$ be
a measurable map of {\it X} into itself preserving the measure $\mu$,
and let $k\ge 3$.
Then
$$
  \liminf_{n\to \infty} \max
  \{ d(T^{n}x,x), d(T^{2n}x,x), \dots, d(T^{(k-1)n}x,x) \} = 0.
$$
for almost all $x\in X$.
}

  A. Behrend  \cite{Be} obtained the following lower bound for  $a_3(N)$
$$
  a_3(N) \gg \exp (-C (\log N)^{\frac{1}{2}} ),
$$
where $C$ is an absolute constant.
A lower bound on $a_k(N)$ for an arbitrary $k$ is given in
\cite{Ra}.

Unfortunately, Szemeredi's methods give very weak upper estimates for
$a_k(N)$.
The ergodic approach gives no estimates at all.
  Only in  2001
W.T. Gowers \cite{Gow_m}
obtained a quantitative result concerning the rate at which $a_k(N)$ approaches zero for $k\ge 4$.
He proved the following theorem.
\\
\Th {\it Let $\delta >0$, $k\ge 4$ and $N\gg \exp \exp (C
\delta^{-K})$, where
$C,K > 0$  is absolute constants.
Let $A\subseteq \{ 1,2,\dots ,N \} $ be a set
of cardinality at least
$\delta N$. Then $A$ contains an arithmetic progression
of length  $k$.}
\\
In other words, W.T. Gowers proved that, for any $k\ge 4$, we have
$a_k(N) \ll 1/ (\log \log N)^{c_k}$, where constant $c_k$ depends
on $k$ only.

  In this paper, we solve the following problem.
  Consider the two--dimensional lattice $[1,N]^2$ with basis  $\{(1,0)$, $(0,1)\}$.
Let
$$
  L(N) = \frac{1}{N^2} \max \{~ |A| ~:~ A\subseteq [1,N]^2 ~\mbox{ and }~
$$
$$
  A \mbox { contains no triples of the form } \{ (k,m),~ (k+d,m),~ (k,m+d)  \}
$$
\begin{equation}\label{tri}
    \mbox { with positive } d \}.
\end{equation}
 A triple from (\ref{tri}) will be called a
 {\it "corner".}
In \cite{Sz2,Fu}, it was proven that  $L(N)$ tends to   $0$ as $N$
tends to infinity.
W.T. Gowers (see \cite{Gow_m}) asked the question
of what is the rate of convergence of $L(N)$ to $0$.

The following theorem was proven in \cite{Shkr_tri,Shkr_tri_DAN} (see also \cite{Sol,Vu}).
\\
\Th {\it Let $\delta>0$ and $N\gg \exp \exp \exp ( \d^{-C} )$, where
$С>0$ is an absolute constant. Let  $A$ be a subset of
$\{1,\dots,N\}^2$ of cardinality at least  $\delta N^2$.
Then $A$ contains a corner.}
\label{main_th_p}
\\
The question on upper estimates for $L(N)$ in the group $\mathbf{F}_3^n$
was considered in \cite{Green_BCC}.
The main result of this paper is the following theorem. \\
\Th {\it Let $\delta>0$, and $N\gg \exp \exp ( \d^{-c} )$, where
$c>0$ is an absolute constant. Let  $A$ be a subset of
$\{1,\dots,N\}^2$ of cardinality at least  $\delta N^2$.
Then $A$ contains a corner.}
\label{main_th}

Thus, we prove the estimate $L(N) \ll 1/ (\log \log N)^{C_1}$,
where $C_1 = 1/c$.

\Note
The constant $c$ in Theorem \ref{main_th} might  be taken as $73$.

The constructions which we use develop the approach of \cite{Bu,Gow_m,Shkr_tri}.

  The proof of Theorem \ref{main_th} is contained in   \S 3,4,5,6
  and  proceeds by an iteration scheme.

Let $A$ be a set, $A\subseteq E_1 \m E_2$, where $E_1, E_2 \subseteq {\bf Z}^2$.
At each step of  our procedure we prove the
following :  either $A$ is "sufficiently regular" or its "density" can be increased.
A suitable definition of  "sufficiently regular"\, sets
(so--called uniform sets) is one of the main aims of our proof.

If $A$ is a random set and $A$ has cardinality $\delta N^2$, then $A$ contains
approximately $\d^3 N^3$ corners. We shall say $A$ is regular
(or in other words {\it $\a$--uniform}) if $A$ contains
the same approximate number of corners.

Let us consider the following example. Let $A$ be a set of the form $E_1\times E_2$,
where $E_1$, $E_2$ are two random sets. Denote by $\beta_i$ the density of the set
$E_i$, and set $\beta_1 \beta_2 = \d$.
Since each $E_i$ has small Fourier coefficients,
so does $A$.
On the other hand, the number of corners in $A$
equals $\beta_1^2 \beta_2^2 N^3 = \d^2 N^3 \neq \d^3 N^3$.
So, if $A$ has small Fourier coefficients, then $A$ might not be regular (uniform).

Let $E_1, E_2$ be subsets of $\L$,  where $\L \subseteq {\bf Z}$ to be chosen later.
Let $A$ be a subset of $E_1 \m E_2$ of cardinality $\d |E_1| |E_2|$.
We shall say that $A$ is {\it rectilinearly $\a$--uniform} if, roughly speaking,
the number of quadruples $\{ (x,y), (x+d,y), (x,y+s), (x+d,y+s) \}$ in $A^4$ is at most
$(\d^4 + \a) |E_1|^2 |E_2|^2$
(in fact we need a slightly different definition of $\a$--uniformity, which depends on the set $\L$).
In \S 3 we prove that if $E_1$, $E_2$ has small Fourier coefficients and $A$ is
rectilinearly $\a$--uniform, then $A$ has about the expected number of corners.

Suppose $A$ fails to be rectilinearly $\a$--uniform.
We shall show in \S 4 that $A$ has increased density
$\d + c(\d)$ on some product set $F_1 \m F_2$, $F_1 \subseteq E_1$, $F_2 \subseteq E_2$.
To obtain this we need  Proposition \ref{na_case_pr}, which was proven by Ben Green in \cite{Green_BCC}.
A similar proposition was proven in \cite{Shkr_tri} with worse bounds.

 Unfortunately, the structure of $F_1 \m F_2$ need not be regular.
 To make it regular, we pass to a subset of $\L$, say, $\L'$ and an integer vector $\v{t}=(t_1,t_2)$ such that
 $(F_1 - t_1) \cap \L'$, $(F_2 - t_2) \cap \L'$ has small Fourier coefficients.

 This is attained by a further iteration procedure.
 Suppose that $F_1 \m F_2$ is not good; then either  $F_1$ or $F_2$ has a large Fourier coefficient.
 This may be used to find a subset of $\L$, say, $\L_1$ such that
 some sort of density (so--called {\it index}, see \S 5)
 of $F_1 \m F_2 $ in $\L_1 \m \L_1$ increases.
 This can only occur  finitely many times.

 We are now in the situation we started with, but $A$ has a larger density
 and we iterate the procedure.
 This also can only occur finitely many times.
 In \S 6 we combine the arguments from the earlier sections and show that they give the bound
 that we stated in Theorem \ref{main_th}.

  The main difference between this paper and \cite{Shkr_tri} consists in the following:
  in \cite{Shkr_tri} we chose $\L$ to be   an {\it arithmetic progression},
  whereas here we put $\L$
  to be
  a so--called {\it Bohr set} (see \cite{Bu, Green_Sz_Ab} and others).
  This choice turns out  to be more economical than dealing with progressions.
  Note that the best upper bound for $a_3 (N)$ was proven by J. Bourgain in \cite{Bu} using
  exactly
  these very sets.
  The properties of Bohr sets will be considered in \S 2.


  At the last section several applications of Theorem \ref{main_th} in the theory of
  dynamical systems will be
  obtained.

  The author is grateful to Professor N.G. Moshchevitin
  for constant attention to this work and to
  Professor Ben Green for helpful conversations
  and
  ideas.

\refstepcounter{section}

\begin{center}
\textsc{ \arabic{section}. On Bohr sets.}
\end{center}

Let $A$ be a subset of ${\bf Z}$.
It is very convenient to write $A(x)$ for such a function.
Thus $A(x)=1$ if $x\in A$ and $A(x)=0$ otherwise.

One of the crucial moments in \cite{Bu} was the notion of Bohr set.

  Let $N$ and $d$ be natural numbers, $\eps>0$ be a real number and $\theta = (\theta_1,\dots, \theta_d) \in {\bf T}^d$.

\Def Define the Bohr set  $\Lambda = \Lambda_{\theta, \eps, N}$ by
$$
  \Lambda_{\theta, \eps, N} = \{ n \in \mathbf{Z} ~|~ |n| \le N,~ \| n \theta_j \| < \eps \mbox{ for } j=1,\dots, d \}
$$

We shall say that the vector $\theta \in {\bf T}^d$ is  {\it generative vector} of Bohr set $\L$.
The number  $d$~ is called  {\it dimension } of Bohr set $\L$ and is denoted by $\dim \L$.
If $M = \L + n$, $n\in {\bf Z}$ is a translation of  $\L$, then, by definition, put $\dim M = \dim \L$.

Another construction of Bohr set (so--called {\it smoothed } Bohr set)
was given in \cite{Tao_lect} and  \cite{Green_Sz_Ab}.

\Def
Let $0<\kappa<1$ be a real number.
A Bohr set $\Lambda = \Lambda_{\theta, \eps, N}$ is called {\it regular}, if
for an arbitrary $\eps'$, $N'$ such that
$$
   | \eps - \eps'|  < \frac{\k}{100 d} \eps \quad \mbox{ and } \quad
   | N - N'|  < \frac{\k}{100 d} N
$$
we have
$$
  1-\k < \frac{|\L_{\theta,\eps',N'}| }{ |\L_{\theta,\eps,N} | } < 1+ \k \,.
$$

We need  several results concerning Bohr sets (see \cite{Bu}).

\Lemma
{\it
Let $\L_{\theta,\eps,N}$ be a Bohr set, $\theta \in {\bf T}^d$.
Then
$$
  |\L_{\theta,\eps,N}| \ge \frac{1}{2} \eps^d N \,.
$$
}
\label{l:Bohr_est}

\Lemma
{\it
Let $0 < \k < 1$ be a real number,
and
$\L_{\theta,\eps,N}$ be a Bohr set.
Then there exists a pair $(\eps_1, N_1)$ such that
$$
  \frac{\eps}{2} < \eps_1 < \eps \quad \mbox{ and } \quad
  \frac{N}{2} < N_1 < N \,,
$$
and such that  $\L_{\theta,\eps_1,N_1}$ is a regular Bohr set.
}
\label{l:Reg_B}

{\it All Bohr sets will be regular in the article.}

\Def Let $f,g$ be functions from $\mathbf{Z}$ to $\bf{C}$.
By $f * g$ define the function
$$
  ( f*g )(n) = \sum_{n\in {\bf Z}} f(s) \overline{ g(n-s) }
$$

\Def
Let $\eps \in (0,1]$ be a real number, and $\L_{\theta,\eps_0,N_0}$ be a Bohr set,
$\theta = (\theta_1, \dots, \theta_d)$.
A regular Bohr set $\L' = \L_{\theta', \eps' , N'}$ is called
{\it $\eps$ attendant }
of $\L$
if
$\theta' = (\theta_1, \dots, \theta_d, \theta_{d+1}, \dots, \theta_{d+k})$, $k\ge 0$,
$\eps \eps_0 /2 \le \eps' \le \eps \eps_0$, $\eps N_0 /2 \le N' \le \eps N_0$.
Lemma \ref{l:Reg_B} implies that for an arbitrary Bohr set there exists its $\eps$ attendant.

We shall consider that $k=0$ unless stated otherwise.

Let $n$ be a natural number, and $\L$ be a Bohr set.
 We shall say that a Bohr set $\L'$ is $\eps$ attendant of $\L+n$,
 if $\L'$ is $\eps$ attendant of $\L$.

The following lemma is also due to J. Bourgain \cite{Bu}.
We give his proof for the sake of completeness.

\Lemma
{\it
Let $\kappa>0$ be a real number, $\theta \in {\bf T}^d$,
$\L = \L_{\theta, \eps, N}$ be a regular Bohr set,
and $\L' = \L_{\theta, \eps', N'}$ its $\kappa / (100 d)$ attendant.
Then the number of $n's$ such that $( \L * \L' ) (n) > 0$ does not exceed  $|\L|(1+\kappa)$,
the number of $n's$ such that $(\L * \L') (n) = |\L'|$
is greater than $|\L|(1-\kappa)$
and
\begin{equation}\label{2k}
  \Big\| \frac{1}{|\L'|} (\L * \L')(n) - \L(n) \Big\|_1 < 2\kappa |\L| \,.
\end{equation}
\label{l:L_pm}
}
\Proof
If $(\L * \L^{'}) (n) > 0$, then there exists $m$ such that
\begin{equation}\label{tM1}
  |m| \le \frac{\k}{100d} N, \quad |n-m| \le N
\end{equation}
and
\begin{equation}\label{tM2}
    \| m \theta_j \| < \frac{\k}{100d} \eps , \quad  \| (n-m) \theta_j \| < \eps, ~ j=1,\dots,d
\end{equation}
Using (\ref{tM1}) and (\ref{tM2}), we get
\begin{equation}\label{}
  |n| \le
            \Big( 1 + \frac{\k}{100d} \Big) N \quad \mbox{ and } \quad
  \| n \theta_j \|
             < \Big( 1 + \frac{\k}{100d} \Big) \eps, ~ j=1,\dots,d
\end{equation}
It follows that
\begin{equation}\label{L+}
  n\in \L^{+} := \L_{\theta,( 1 + \frac{\k}{100d} ) \eps, ( 1 + \frac{\k}{100d} ) N}  \,.
\end{equation}
By Lemma \ref{l:Reg_B} we have $|\L^{+}| \le (1+\k) |\L|$.

On the other hand, if
\begin{equation}\label{L-}
  n\in \L^{-} := \L_{\theta,( 1 - \frac{\k}{100d} ) \eps, ( 1 - \frac{\k}{100d} ) N} \,,
\end{equation}
then $(\L * \L^{'}) (n) = |\L'|$.
Using Lemma \ref{l:Reg_B}, we obtain $|\L^{-}| \ge (1-\k) |\L|$.

Let us prove (\ref{2k}).
We have
$$
  \Big\| \frac{1}{|\L'|} (\L * \L')(n) - \L(n) \Big\|_1
    =
  \Big\| \frac{1}{|\L'|} (\L * \L')(n) - \L(n) \Big\|_{ l^1 ( \L^{+} \setminus \L^{-}) }
$$
$$
  \le
  |\L^{+}| - |\L^{-}| < 2 \k |\L|
$$
as required.

\Cor Lemma  \ref{l:L_pm} implies that  $|\L| \le |\L + \L'| \le (1+2\k) |\L|$.

\Note Let $\L^x (n) = \L (n-x)$.
Since  $(\L^x * \L') (n) = (\L * \L')(n-x)$, it follows that (\ref{2k}) takes place for translations   $\L + x$.

\Def
By $\L^{+}$ and $\L^{-}$ denote the Bohr sets defined in (\ref{L+}) and (\ref{L-}), respectively,
$\L^{-} \subseteq \L \subseteq \L^{+}$.
By Lemma \ref{l:L_pm} we have $|\L^{+}| \le |\L|(1+\kappa)$
and $|\L^{-}| \ge |\L|(1-\kappa)$.
Note that for any $s\in \L^{'}$, we get $\L^{-} \subseteq \L + s$.

Suppose  $\L \subseteq {\bf Z}$ is  a Bohr set, and $\v{x} = (x_1, x_2)$ belongs to ${\bf Z^2}$.
By $\L + \v{x}$ denote the set $(\L+x_1) \m (\L+x_2) \subseteq {\bf Z^2}$.
Let $\v{n} \in {\bf Z^2}$.
Let $\L(\v{n})$ denote the characteristic function of $\L \m \L$.
We shall write  $\v{s} \in \L$, $\v{s} = (s_1, s_2)$, if $s_1 \in \L$ and $s_2 \in \L$.

\Lemma
{\it
  Suppose $\L$ is a Bohr set, $\L^{'}$ is its $\eps$ attendant, $\eps = \k/ (100 d)$, $\v{x}$ is a vector,
  and $E \subseteq {\bf Z^2}$.
  Then
  \begin{equation}\label{}
    \Big| \d_{\L+ \v{x}} (E) -  \frac{1}{|\L|^2} \sum_{\v{n} \in \L+ \v{x}} \d_{\L^{'} + \v{n}} (E)  \Big|
        \le 4\k \,.
  \end{equation}
}
\label{l:smooth_d}
\Proof
We have
$$
  \sigma =
  \frac{1}{|\L|^2} \sum_{\v{n} \in \L+ \v{x}} \d_{\L^{'} + \v{n}} (E)
    =
  \frac{1}{|\L|^2 |\L^{'}|^2} \sum_{\v{s}} E(\v{s}) \sum_{\v{n}} \L(\v{n} - \v{x}) \L^{'} (\v{s} - \v{n})
$$
$$
    =
  \frac{1}{|\L|^2 |\L^{'}|^2} \sum_{\v{s}} E(\v{s}) \sum_{\v{n}} \L(\v{n} ) \L^{'} (\v{s} - \v{x} - \v{n})
$$
Using Lemma \ref{l:L_pm}, we get
$$
  \sigma = \frac{1}{|\L|^2} \sum_{\v{s}} E(\v{s}) \L( \v{s} - \v{x} ) + 4 \vartheta \k
         = \d_{\L+ \v{x}} (E) + 4 \vartheta \k \,,
$$
where  $|\vartheta | \le 1$.
This completes the proof.

\Note
Clearly, the one--dimension analog of Lemma \ref{l:smooth_d} takes place.

Let $\L_1 = \L_{\theta_1,\eps_1,N_1}$, $\L_2 = \L_{\theta_2,\eps_2,N_2}$ be two Bohr sets.
We shall write $\L_1 \le \L_2$, if $\theta_1 = \theta_2$, $\eps_1 \le \eps_2$ and $N_1 \le N_2$.

Note that if $\L_1 \le \L_2$, then an arbitrary $\eps$ attendant of $\L_1$ is $\eps$ attedant of $\L_2$.

\refstepcounter{section}

\begin{center}
\textsc{ \arabic{section}. On $\alpha$--uniformity.}
\end{center}

%
%
%
%

Let $f$ be a function from $\mathbf{Z}$ to  $\mathbf{C}$.
By  $\F{f}(x)$ denote the Fourier transformation of $f$
$$
  \F{f}(x) =  \sum_{s\in {\bf Z}} f(s) e( -sx) ,
$$
where
$e(x) = e^{2\pi ix}$.
\\
We shall use the following basic facts
\begin{equation}\label{F_Par}
    \sum_{s\in {\bf Z}} |f(s)|^2 = \int_0^1 |\widehat{f} (x)|^2 dx
\end{equation}
\begin{equation}\label{F_Par_sc}
    \sum_{s\in {\bf Z}} f(s) \overline{g(s)}= \int_0^1 \F{f}(x) \overline{\F{g}(x)} dx
\end{equation}
\begin{equation}\label{svertka}
    \sum_{k\in {\bf Z}} | \sum_{s\in {\bf Z}} f(s) \overline{g(s-k)} |^2
        = \int_0^1 |\widehat{f} (x)|^2 |\widehat{g} (x)|^2 dx
\end{equation}

  Let $\Lambda$ be a Bohr  set,  and $A$ be an arbitrary subset of $\Lambda$.
  Let $|A| = \d |\Lambda|$.
  Define the {\it balanced} function of $A$ to be $f(s) = ( A(s) - \d) \L (s) = A(s) - \d \L (s)$.
\\
Let ${\bf D}$ denote the closed disk of radius $1$ centered at $0$
in the complex plane. Let $R$ be an arbitrary set.
We write  $f: R \to {\bf D}$ if $f$ is zero outside $R$.
\\
The following definition is due to Gowers \cite{Gow_m}.
\\
\Def
  A function $f : \L \to {\bf D}$ is called $\a$--uniform if
\begin{equation}
 \| \F{f} \|_{\infty}  \le \a |\Lambda|
\label{a_d0}
\end{equation}
\label{d:basic}

We say that  $A$ is $\a$--uniform if its balanced function  is.

We shall write $\int$ instead of $\int_0^{1}$ and $\sum_s$ instead of $\sum_{s\in {\bf Z}}$.

Let us prove an analog of Lemma 2.2 from \cite{Gow_m}.

\Lemma {\it Let $\L$ be a Bohr set, and let  $f : \L \to {\bf D}$ be $\a$--uniform function.
Then we have
$$
   \sum_k |\sum_s f(s) \overline{g(s-k)} |^2
     \le \a^2 |\L|^2 \| g \|^2_2,
$$
for an arbitrary function $g$, $g: {\bf Z} \to {\bf D}$.
\label{aux_l}
}
\\
\Proof
By  (\ref{svertka}) we get
\begin{equation}\label{s:1}
  \sum_k |\sum_s f(s) \overline{g(s-k)} |^2 =
  \int |\widehat{f} (x)|^2 |\widehat{g} (x)|^2 dx.
\end{equation}
Since the function $f$ is $\a$--uniform, it follows that $\| \F{f} \|_{\infty}  \le \a |\Lambda|$.
Using this inequality and (\ref{F_Par_sc}), we have
\begin{equation}\label{}
   \sum_k |\sum_s f(s) \overline{g(s-k)} |^2 \le \a^2 |\L|^2 \int |\widehat{g} (x)|^2 dx
   \le \a^2 |\L|^2 \| g \|^2_2
\end{equation}
This completes the proof.

\Cor  {\it Let $S$  be a set, and $\L^{'}$ be a Bohr set.
Suppose $E\subseteq \L^{'}$ is $\a$--uniform, and $E$ have the cardinality $\d |\L^{'}|$.
Let  $g$ be a function from $S$ to ${\bf D}$.
Then for all but $\a^{2/3} |S|$ choices of $k$ we have
$$
   \Big| (E * g) (k) - \d ( \L^{'} * g)(k) \Big| \le \a^{2/3} |\L^{'}| \,.
$$
\label{c:sv}
}
\\
Let  $f$ be the balanced function of $E\cap \L^{'}$.
Using Lemma \ref{aux_l}, we get
\begin{equation}
     \sum_k | ( E * g) (k) - \d ( \L^{'} * g)(k) |^2 =
     \sum_k |\sum_s f(s) \overline{g(s-k)} |^2 \le
\end{equation}
\begin{equation}\label{}
     \le \a^2 |\L^{'}|^2 \| g \|^2_2 \le \a^2 |\L^{'}|^2 |S| \,.
\end{equation}
This concludes the proof.

  By $\v{e}_1$ and $\v{e}_2$ define the vectors $(1,0)$ and $(0,-1)$.

  Let $\L_1$ and $\L_2$ be Bohr sets,
  and  $E_1\m E_2$ be  a subset of $\L_1 \m \L_2$.
Suppose  $f:\L_1\m \L_2 \to {\bf D}$ is a function.

\Def  Let $\alpha$ be a real number, $\a \in [0,1]$.
A function $f : E_1 \m E_2\to {\bf D}$ is called {\it rectilinearly  $\a$--uniform} if
\begin{equation}
  \sum_{\vec{s},u} \sum_r f(\v{s}) \overline{f(\v{s}+u\v{e}_2)}
  \overline{f(\v{s}+r\v{e}_1)} f(\v{s} + u\v{e}_2+ r\v{e}_1) \le \a |E_1|^2 |E_2|^2.
\label{a_d1}
\end{equation}
\label{a_d1_ref}
Let $f(k,m) = f(k\v{e}_1 + m\v{e}_2)$.
Note that the function $f$ is $\a$--uniform iff
\begin{equation}
  \sum_{m,p} |\sum_k f(k,m) \overline{f(k,p)}|^2 \le \a |E_1|^2 |E_2|^2.
\label{a_d2}
\end{equation}

  Let $A$ be a subset of  $E_1\m E_2$, $|A| = \d |E_1| |E_2|$.
Define the {\it balanced} function of $A$ to be $f(\vec{s})=( A(\vec{s})- \delta ) \cdot (E_1\m E_2) (\v{s})$.
  We say that $A\subseteq E_1\m E_2$ is rectilinearly $\a$--uniform if its balanced function is.

  Let $f$ be an arbitrary function, $f : {\bf Z}^2 \to \mathbf{C}$.
Define $\| f \|$ by the formula
\begin{equation}
  \| f \| =
  | \sum_{\vec{s},u} \sum_r f(\v{s})
  \overline{f(\v{s}+u\v{e}_2)}
  \overline{f(\v{s}+r\v{e}_1)} f(\v{s} + u\v{e}_2+ r\v{e}_1) |^{\frac{1}{4}}
\label{norm}
\end{equation}

\Lemma {\it $\| \cdot \|$ is a norm.}
\label{l:norm}
\\
\Proof See \cite{Shkr_tri}.

%
%
%
%
%

\Def
Let $\L$ be a Bohr set, $Q \subseteq \L$, $|Q| = \d |\L|$,
$\a,\eps$ are positive numbers,
and $\L'$ be  $\eps$  attendant set of $\L$.
Consider the set
$$
   B = \{ m \in \L  ~|~  \| (Q \cap (\L'+m) - \d (\L'+m) ) \F{} ~ \|_{\infty}  \ge \a |\L'| \} \,.
$$
A set $Q$ is called {\it $(\a,\eps)$--uniform} if
\begin{equation}\label{c:sm_B}
  |B| \le \a |\L| \,,
\end{equation}
\begin{equation}\label{c:sm_sq}
  \frac{1}{|\L|} \sum_{m\in \L} |\d_{\L'+m} (Q) - \d |^2 \le \a^2 \,.
\end{equation}
and
\begin{equation}\label{c:sm_total_F}
  \| ( Q \cap \L - \d \L ) \F{} ~ \|_{\infty} \le \a |\L| \,.
\end{equation}
Certainly, this definition depends on $\L$ and $\L'$.
We do not assume that $\L'$ has the same generative vector as $\L$.

\Note
Let
$$
  B^* = \{ m \in \L  ~|~  | \d_{\L'+m} (Q) - \d |  \ge \a^{2/3}  \} \,.
$$
Condition (\ref{c:sm_sq}) implies  that $|B^*| \le \a^{2/3} |\L|$.

\Note
Condition (\ref{c:sm_total_F}) is not so important as (\ref{c:sm_B}) and (\ref{c:sm_sq}).
The inequality
$$
  \| ( Q \cap \L - \d \L ) \F{} ~ \|_{\infty} \le 4 \a |\L|
$$
follows from (\ref{c:sm_B}), (\ref{c:sm_sq}) (see Statement \ref{l:l2}).

Let $\L_1$, $\L_2$ be Bohr sets, $\L_1 \le \L_2$,
$\eps>0$ be a real number, and $\L'$ be $\eps$ attendant of $\L_1$.
Let also $E_1$, $E_2$ be subsets of $\L_1$, $\L_2$, respectively,  and
$|E_1| = \beta_1 |\L_1|$, $E_2 = \beta_2 |\L_2|$.

\Def
A function $f : \L_1 \m \L_2 \to {\bf D} $ is called {\it rectilinearly $(\a, \eps)$--uniform}
if
$$
  \| f \|_{\L_1\m \L_2, \eps}^4 = \sum_{i\in \L_1} \sum_{j\in \L_2}
                                  \sum_k \sum_{m,u}
                                  \L' (m-k-i) \L' (u-k-i) ~\m
$$
\begin{equation}\label{f:f_for_norm}
                                  |\sum_r \L' (k+r-j) f(r,m) f(r,u) |^2
                                  \le
                                  \a \beta_1^2 \beta_2^2 |\L'|^4 |\L_1|^2 |\L_2| \,.
\end{equation}


Let $\L_1$, $\L_2$ be Bohr sets, $\L_1 \le \L_2$, and
$\L'$ be $\eps$ attendant of $\L_1$.
Suppose that $\L'_\eps$ is $\eps$ attendant of $\L'$.
Let also $E_1$, $E_2$ be subsets of $\L_1$, $\L_2$, respectively,  and
$|E_1| = \beta_1 |\L_1|$, $E_2 = \beta_2 |\L_2|$.

\Def Let  $A \subseteq E_1 \m E_2$, $|A| = \d \beta_1 \beta_2 |\L_1|  |\L_2|$,
and
$f(\v{s}) = A(\v{s}) - \delta (E_1 \m E_2) (\v{s})$.
Let $f_l (\v{s}) = f(s_1 + l, s_2) \L'(s_1) $, $l\in \L_1$.                                              
Consider the set
$$
  B = \{ l \in \L_1 ~|~ \| f_l \|^4_{\L' \m \L_2,\eps} >
         \a \beta_1^2 \beta_2^2 |\L'_\eps|^4 |\L'|^2 |\L_2| \} \,.
$$
$A$ is called {\it rectilinearly $(\a, \a_1, \eps)$--uniform}
if $|B| \le \a_1 |\L_1|$.
\\
Note that
$$
  \| f_l \|^4_{\L' \m \L_2,\eps} = \sum_{i\in \L'} \sum_{j\in \L_2}
                                   \sum_k \sum_{m,u}
                                   \L'' (m-k-i) \L'' (u-k-i)
                                   |\sum_r \L'' (k+r-j) f_l (r,m) f_l(r,u) |^2
$$
$$
                                  =
                                   \sum_{i\in \L' + l} \sum_{j\in \L_2}
                                   \sum_k \sum_{m,u}
                                   \L'' (m-k-i) \L'' (u-k-i) ~\m
                                   |\sum_r \L'' (k+r-j) \t{f}_l (r,m) \t{f}_l(r,u) |^2 \,,
$$
where $\L'' = \L' (\eps)$ and $\t{f}$ is a restriction of  $f$ to $(\L'+l) \m \L_2$.

\label{def:rectangle}
\Note
We need parameter $\a_1$ to decrease the constant $c$ in Theorem \ref{main_th}.
To obtain Theorem \ref{main_th} with $c$ equals, say, $1000$, one can put  $\a_1 = \a$.

\Lemma
{\it
  Let $\L$ be a Bohr set.
  Suppose $\L'$ is  $\eps$ attendant of $\L$,
  $\L''$ is $\eps$ attendant $\L'$  and $\eps^2$ attendant of $\L$,
  $\eps = \a^2 / 4 (100 d)$,
  $Q \subseteq \L$, $|Q| = \d \L$, and $\a>0$.
  Let
  $$
    \Omega_1 = \{ s\in \L ~|~ | \d_{\L'+s} (Q) - \d | \ge 4 \a^{1/2}
    \mbox{ or }~ \frac{1}{|\L'|} \sum_{n\in \L'+s} | \d_{\L'' + n} (Q) - \d |^2 \ge 4 \a^{1/2} \} \,.
  $$
  $$
    \Omega_2 = \{ s\in \L ~|~ \| (Q \cap (\L' + s)  - \d (\L'+s) ) \F{} ~ \|_{\infty} \ge 4 \a^{1/4} |\L'| \}  \,.
  $$
1)
  If
  \begin{equation}\label{c:T1}
    \frac{1}{|\L|} \sum_{n\in \L} | \d_{\L'' + n} (Q) - \d |^2 \le \a^2 \,,
  \end{equation}
  then   $ |\Omega_1| \le 4 \a^{1/2} |\L|$.
  \\
2)
  If
  \begin{equation}\label{c:T2}
    \Omega^* = \{ s\in \L ~|~ \| (Q \cap (\L'' + s) - \d (\L''+s) ) \F{} ~ \|_{\infty} \ge \a |\L''| \}
  \end{equation}
  has the cardinality at most $\a |\L|$,
  then   $ |\Omega_2| \le 4 \a^{1/2} |\L|$.
  \\
3)
  Suppose $Q$ is $(\a,\eps^2)$--uniform subset of $\L$.
  Let
  $$
    \t{\Omega} = \{ s\in \L ~|~ \mbox{ Set } ( Q - s ) \cap \L' \mbox{ is not }
                                (8 \a^{1/4},\eps)\mbox{--uniform } \} \,.
  $$
  Then $|\t{\Omega}| \le 8 \a^{1/2} |\L|$.
}
\label{l:intermediate}
\\
\Proof
Let us prove $1)$.
Let $\d'_n = \d_{\L' + n} (Q)$, $\d''_n = \d_{\L'' + n} (Q)$, $\k = \a^2 / 4$, and $\epsilon = \a^{1/2}$.
Сonsider the sets
$$
  B_s = \{  n \in \L' + s ~|~ | \d''_n - \d | \ge \epsilon \} \,,
  G_s = \{  n \in \L' + s ~|~ | \d''_n - \d | < \epsilon \} \,, s\in \L
$$
and sets
$$
  B = \{ s \in \L ~|~ |B_s| \ge \epsilon |\L'| \} \,,
  G = \{ s \in \L ~|~ |B_s| < \epsilon |\L'| \}
$$
If $s\in G$, then $|B_s| < \epsilon |\L'|$.
Using Lemma \ref{l:smooth_d}, we have
$$
  | \d'_s - \d | \le \Big| \frac{1}{|\L'|} \sum_{x\in \L'+s} \d''_x   -   \d \Big| + 4 \k
  \le
    \frac{1}{|\L'|}  \sum_{x\in \L'+s} |\d''_x   -   \d | + 4 \k
  \le
$$
\begin{equation}\label{tmp:17:17}
  \le
    \frac{1}{|\L'|}  \sum_{x\in B_s} |\d''_x   -   \d | + \frac{1}{|\L'|}  \sum_{x\in G_s} |\d''_x   -   \d | + 4 \k
  <
      \epsilon + \frac{\epsilon |G_s|}{|\L'|} + 4 \k
  \le
      4 \epsilon \,.
\end{equation}
Besides that for $s\in G$, we get
\begin{equation}\label{tmp:17:17+}
  \frac{1}{|\L'|} \sum_{x\in \L'+s} |\d''_x   -   \d |^2
    \le
        \frac{1}{|\L'|} \sum_{x\in B_s} |\d''_x   -   \d |^2
            +
        \frac{1}{|\L'|} \sum_{x\in G_s} |\d''_x   -   \d |^2
            \le
            \epsilon + \epsilon^2 \le 2 \epsilon \,.
\end{equation}
Let us estimate the cardinality of $B$.
We have
$$
 \a^2 \ge
        \frac{1}{|\L|} \sum_{s\in B} |\d''_s - \d|^2
    \ge
        \frac{1}{|\L'| |\L|} \sum_{s\in B} \sum_{n \in \L'+s} |\d''_n - \d|^2 - 4\k
    \ge
$$
$$
    \ge
        \frac{1}{|\L'| |\L|} \sum_{s\in B} \sum_{n \in B_s} |\d''_n - \d|^2 - 4 \k
    \ge
        \frac{|B| \epsilon^3 |\L'|}{|\L'| |\L|}  - 4\k \,.
$$
It follows that, $|B| \le 4 \a^{1/2} |\L|$.
Using  (\ref{tmp:17:17}), (\ref{tmp:17:17+}) we get
$\Omega_1 \subseteq B$ and  $1)$ is proven.
To prove $2)$ it suffices to note that                                                                  
$$
  \frac{1}{|\L||\L'|} \sum_{s\in \L} \| (Q \cap (\L'' + s)  - \d (\L''+s) ) \F{} ~  \|_{\infty}
      =
  \frac{1}{|\L||\L'|} \sum_{s\in \Omega^*} \| (Q \cap (\L'' + s) - \d (\L''+s) ) \F{} ~  \|_{\infty}
    +
$$
$$
    +
  \frac{1}{|\L||\L'|} \sum_{s\in ( \L \setminus \Omega^* ) } \| (Q \cap (\L'' + s) - \d (\L'' + s)) \F{} ~  \|_{\infty}
      \le
            \a + \frac{\a |\L'|}{|\L||\L'|} |\L \setminus \Omega^*|
      \le 2\a \,.
$$
and define the sets $B'_s$, $G'_s$, $B'$, $G'$ ~:~
$$
  B'_s = \{  n \in \L' + s ~|~ \| (Q \cap (\L'' + s) - \d (\L''+s) ) \F{} ~ \|_{\infty} \ge \epsilon_1 |\L''| \} \,,
$$
$$
  G'_s = \{  n \in \L' + s ~|~ \| (Q \cap (\L'' + s) - \d (\L''+s) ) \F{} ~ \|_{\infty} < \epsilon_1 |\L''| \}, \,~ s\in \L \,.
$$
$$
  B' = \{ s \in \L ~|~ |B_s| \ge \epsilon_1 |\L'| \}
  \quad \mbox{ and } \quad
  G' = \{ s \in \L ~|~ |B_s| < \epsilon_1 |\L'| \} \,,
$$
where
$\epsilon_1 = \a^{1/4}$.
After that we can apply the same arguments as above, using
Lemma \ref{l:L_pm} instead of Lemma \ref{l:smooth_d}.

Let us prove $3)$.
Since  $Q$ is $(\a,\eps^2)$--uniform subset of $\L$, it follows that $Q$ satisfies (\ref{c:T1}).
Also we have $|\Omega^*| \le \a |\L|$ and $|B|, |B'| \le 4 \a^{1/2} |\L|$.
It is easily shown that for all $s \notin B \cup B'$ the set $( Q - s ) \cap \L'$ is
$(8 \a^{1/4},\eps)$--uniform.
This completes the proof.

In the same way we can prove

\St
{\it
  Let $\L$ be a Bohr set,  and $E\subseteq \L$, $|Q| = \d |\L|$ be  $(\a, \eps)$--uniform,
  $\eps = \a / 4 (100 d)$.
  Then
  \begin{equation}\label{}
   \| (Q \cap \L - \d \L) \F{} ~ \|_{\infty} < 4 \a |\L| \,.
  \end{equation}
}
\label{l:l2}
We will not, however, use this fact.
\\

Let $\L_1$, $\L_2$ be Bohr sets, $\L_1 \le \L_2$, and $E_1 \subseteq \L_1$, $E_2 \subseteq \L_2$,
$|E_1| = \beta_1 |\L_1|$, $|E_2| = \beta_2 |E_2|$.
By  $\mathcal{P}$ denote the $E_1 \m E_2$.
Let $H, G $ be subsets of $\mathcal{P}$.

\Th
{\it
Let $f : \mathcal{P} \to {\bf D}$  be a function.
Suppose that $f$ is rectilinearly
$(\a,\eps)$--uniform,
and
the
sets $E_1$, $E_2$ are
$(\a_0,\eps)$--uniform,
$\a_0 = 2^{-50} \a^2 \beta_1^{12} \beta_2^{12}$,
$\eps = 2^{-10} \eps_0^2$,
$\eps_0 = (2^{-10} \a_0^2) /(100 d)$.
Let also
$\L_1$ be $\eps_0$ attendant of $\L_2$.
Then
\begin{equation}\label{f:f_residual}
  | \sum_{\v{s}\in {\bf Z^2 } } \sum_{r \in {\bf Z} }
  H (\v{s}) G ( \v{s} + r\v{e} ) f(\v{s} + r\v{e}_2) |
  \le 2^5 \a^{1/4} \beta_1^2 \beta_2^2 |\L_1|^2 |\L_2| \,.
\end{equation}
}
\label{t:tmp_tha}
\Proof
Let $\v{e}= \v{e}_1 + \v{e}_2$, $\v{s}=k\v{e}_1 +
m\v{e}_2$, and $\L'$ be $\eps$ attendant of $\L_1$.

Let
$$
   \Omega^{(1)}_1 = \{ s \in \L_1  ~|~ \| (E_1 \cap (\L'+s) - \d (\L'+s) ) \F{} ~ \|_{\infty} \ge \a_0  \} \,, \quad
$$
$$
   \Omega^{(1)}_2 = \{ s \in \L_1  ~|~ | \d_{\L'+s}( E_1) - \beta_1 | \ge \a_0^{2/3} \} \,,
$$
and
$$
  \Omega^{(2)}_1 = \{ s \in \L_2  ~|~ \| (E_2 \cap (\L'+s) - \d (\L'+s) ) \F{} ~ \|_{\infty} \ge \a_0  \} \,, \quad
$$
$$
  \Omega^{(2)}_2 = \{ s \in \L_2  ~|~ | \d_{\L'+s}( E_2) - \beta_2 | \ge \a_0^{2/3} \} \,,
$$
Let also
$\Omega_1 = \Omega^{(1)}_1 \cup \Omega^{(1)}_2$,
and
$\Omega_2 = \Omega^{(2)}_1 \cup \Omega^{(2)}_2$.
By assumption  the sets $E_1$, $E_2$ are $(\a_0, \eps)$--uniform.
It follows that $|\Omega^{(1)}_l| \le \a_0^{2/3} |\L_1| $, $|\Omega^{(2)}_l| \le \a_0^{2/3} |\L_2|$, $l=1,2$.
Hence, $|\Omega_1| \le 2 \a_0^{2/3} |\L_1|$ and $|\Omega_2| \le 2 \a_0^{2/3} |\L_2|$.

Let $g_i(\v{s}) = g_i (k,m) = G (k,m) \L' (k-i) $, $i\in \L_1$, and $h_j(\v{s}) = h_j (k,m) = H (k,m) \L' (m-j)$, $j\in \L_2$.
We have $k \in \L_1$, $m \in \L_2$ and $k+r \in \L_1$ in (\ref{f:f_residual}).                              
It follows that the sum (\ref{f:f_residual}) does not exceed $|\L_1|^2 |\L_2|$.
Let also $\l_i = \L' + i$, and $\mu_j = \L' + j$.
Using Lemma \ref{l:L_pm}, we get
$$
  \sigma_0 = \sum_{\v{s}\in {\bf Z^2 } } \sum_{r \in {\bf Z} }
             H (\v{s}) G ( \v{s} + r\v{e} ) f(\v{s} + r\v{e}_2)
           =
$$
$$
           =
             \sum_{k,m} \sum_{r}
             H (k,m) G ( k + r, m+r ) f(k,m+r) \L_1(k+r) \L_2(m)
           =
$$
$$
             \frac{1}{|\L'|^2}
             \sum_{k,m} \sum_{r}
             H (k,m) G ( k + r, m+r ) f(k,m+r) (\L_1 * \L') (k+r) (\L_2 * \L') (m) + 16 \vartheta_0 \k |\L_1|^2 |\L_2|
           =
$$
\begin{equation}\label{f:begG}
             \frac{1}{|\L'|^2}
             \sum_{i\in \L_1} \sum_{j\in \L_2}
             \sum_{k,m} \sum_r
             h_j (k,m) g_i ( k + r, m+r ) f(k,m+r)
             + 16 \vartheta_0 \k |\L_1|^2 |\L_2| \,,
\end{equation}
where $ |\vartheta_0| \le 1$ and $\k \le 2^{-10} \a_0^2$.
Split the sum $\sigma_0$ as
\begin{equation}\label{f:4+R}
  \sigma_0 = \widetilde{\sigma}_0 + \sigma_0^{'} + \sigma_0^{''} + \sigma_0^{'''} + R\,,
\end{equation}
The sum $\widetilde{\sigma}_0$ is taken over $i \notin \Omega_1, j\notin \Omega_2$,
the sum $\sigma_{0}^{'}$ is taken over $i \in \Omega_1, j\notin \Omega_2$,
the sum $\sigma_{0}^{''}$ is taken over $i \notin \Omega_1, j\in \Omega_2$,
the sum $\sigma_{0}^{'''}$ is taken over $i \in \Omega_1, j\in \Omega_2$
and
$|R| \le 16 \eps |\L_1|^2 |\L_2|$.
Let us estimate  $\sigma_0^{'}$, $\sigma_0^{''}$ and  $\sigma_0^{'''}$.
Rewrite $\sigma_0$ as
\begin{equation}\label{f:f_residual1}
  \sigma_0 = \frac{1}{|\L'|^2} \sum_{i \in \L_1} \sum_{j\in \L_2} \sum_{k,m} \sum_r h_j(k-r,m) g_i(k,m+r) f(k-r,m+r) + R\,.
\end{equation}
Let $i$ and $j$ in the sum (\ref{f:f_residual1}) be fixed.
We have $k\in \l_i$ and $m\in \mu_j$.
Further if  $f(k-r,m+r)$ is not zero, then $k-r \in \L_1$.
It follows that $r\in \l_i - \L_1 = \L' - \L_1 + i$.
The set $\L'$ is $\eps$ attendant of $\L_1$.
Using Lemma \ref{l:L_pm}, we obtain that $r$ belongs to a set of cardinality at most  $2 |\L_1|$.
Hence
\begin{equation}\label{e:sigma'}
  |\sigma_0^{'}|  \le \frac{1}{|\L'|^2} 2 |\Omega_1| \cdot |\L_2| \cdot |\L'|^2 |\L_1|
                  \le 2 \a_0^{2/3} |\L_1|^2 |\L_2| \,.
\end{equation}
In the same way $|\sigma_0^{''}| \le 2 \a_0^{2/3} |\L_1|^2 |\L_2| $ and
$|\sigma_0^{'''}| \le 2 \a_0^{2/3} |\L_1|^2 |\L_2| $.

Take $i$ and $j$ such that   $i\notin \Omega_1$, $j\notin \Omega_2$.
Let  $g (\v{s}) = g_i (\v{s})$, $h (\v{s}) = h_j (\v{s})$, and $\L_1 \m \mu_j = \L_1^{(1)} \m \L_2^{(1)}$,
$\l_i \m \L_2 = \L_1^{(2)} \m \L_2^{(2)}$.
Let
$E_2^{(1)} = E_2 \cap \L_2^{(1)}$, $E_1^{(2)} = E_1 \cap \L_1^{(2)}$,
$\beta_2^{(1)} = |E_2^{(1)}| / |\L_2^{(1)}|$, and $ \beta_1^{(2)} = |E_1^{(2)}| / |\L_1^{(2)}|$.
We have
\begin{equation}\label{}
 \sigma = \sigma_{i,j}
        = \sum_{\v{s}\in {\bf Z^2 } } \sum_{r \in {\bf Z} }
          h (\v{s}) g (\v{s} + r(\v{e}_1 + \v{e}_2)) f(\v{s} + r\v{e}_2) =
 \end{equation}
\begin{equation}\label{}
        = \sum_{k,m} h(k,m) E_2^{(1)} (m) \sum_{r} g(k+r,m+r) f(k,m+r)
\end{equation}
Using  the Cauchy--Bounyakovskiy inequality, we obtain
\begin{equation}\label{}
 |\sigma|^2 \le \| h \|_2^2 \sum_{k,m} E_2^{(1)} (m) |\sum_r g(k+r,m+r) f(k,m+r) |^2 =
\end{equation}
$$
 = \| h \|_2^2 \sum_{k,m} E_2^{(1)} (m) \sum_{r,p} g(k+r,m+r) f(k,m+r) g(k+p,m+p) f(k,m+p) =
$$
$$
 = \| h \|_2^2 \sum_{k,m} E_2^{(1)} (m-r) \sum_{r,u} g(k,m) f(k-r,m) g(k+u,m+u) f(k-r,m+u) =
$$
$$
 = \| h \|_2^2 \sum_{k,m} \sum_u g(k,m) g(k+u,m+u) \sum_r E_2^{(1)} (m-r) f(k-r,m) f(k-r,m+u) =
$$
$$
 = \| h \|_2^2 \sum_{k,m} \sum_u g(k,m) g(k+u,m+u) E_1^{(2)} (k) E_1^{(2)} (k+u) 
$$
\begin{equation}\label{}
  \cdot \sum_r  E_2^{(1)} (m-r) f(k-r,m) f(k-r,m+u)
\end{equation}
We have $k\in \L_1^{(2)}$ and $k-r \in \L_1$.
It follows that $r\in k-\L_1 \in \L_1^{(2)} - \L_1$.
Since  $m-r \in \L_2^{(1)}$, it follows that $m\in \L_2^{(1)} + r \in \L_2^{(1)} + \L_1^{(2)} - \L_1$.
On the other hand $k+u \in \L_1^{(2)}$.
Hence $u\in \L_1^{(2)} - \L_1^{(2)}$ and
$m+u \in \L_2^{(1)} + \L_1^{(2)} - \L_1 + \L_1^{(2)} - \L_1^{(2)}$.
Let $\tilde{\L}_i = \L' + \L' + \L' + \L' + \L_1 + i$.
Then $m, m+u \in \tilde{\L}_i + j = Q_{ij} = Q$.
Using Lemma \ref{l:L_pm} for Bohr set $\L_1$ and its $\eps$ attendant  $\L'$, we obtain
that the cardinality of $\tilde{\L}_i$ does not exceed $5|\L_1|$.
Using the Cauchy--Bounyakovskiy inequality, we get
\begin{equation}\label{}
  |\sigma|^4 \le \| h \|_2^4 \Big( \sum_{k} \sum_{m,u} g(k,m) g(k+u,m+u) \Big)
\end{equation}
$$
                           \cdot
                           \Big(
                                 \sum_{k,m,u} E_1^{(2)} (k) E_1^{(2)} (k+u)
                                 \sum_{r,r' } E_2^{(1)} (m-r) E_2^{(1)} (m-r') \, \times
$$
$$
                                 \times \,
                                              f(k-r,m) f(k-r,m+u) f(k-r',m) f(k-r',m+u)
                           \Big)
$$
\\

Let us estimate  $\sigma^{*} = \sigma^{*}_{ij} = \sum_{k,m,u} g(k,m) g(k+u,m+u)$.
Let $\tilde{E}_2^{(2)} = E_2 \cap Q$.
We have
$$
  \sigma^{*} = \sum_{k,m,u} g(k,m) g(k+u,m+u)
               \le \sum_{k,m,u} E_1^{(2)} (k) E_1^{(2)} (k+u) \tilde{E}_2^{(2)} (m) \tilde{E}_2^{(2)} (m+u)
$$
\begin{equation}\label{}
               = \sum_{k,m,u} E_1^{(2)} (k) E_1^{(2)} (u) \tilde{E}_2^{(2)} (m) \tilde{E}_2^{(2)} (m+u-k) =
\end{equation}
\begin{equation}\label{}
               = \sum_{k,m,u} E_1^{(2)} (k) \tilde{E}_2^{(2)} (m+k) E_1^{(2)} (u) \tilde{E}_2^{(2)} (m+u) =
\end{equation}
\begin{equation}\label{tmp_stupid_est}
               = \sum_m ( E_1^{(2)} * \tilde{E}_2^{(2)} )^2 (m) \,.
\end{equation}
Recall that  $|Q_{ij}| \le 5|\L_1|$.
Lemma \ref{l:L_pm} implies that $m$ in the sum (\ref{tmp_stupid_est})
belongs to a set of cardinality at most $8|\L_1|$.
The expression  (\ref{tmp_stupid_est}) implies that for all $i,j$ we have
\begin{equation}\label{}
  |\sigma^{*}_{ij} | \le 8 |\L'|^2 |\L_1| \,.
\end{equation}
We need a stronger upper bound for $\sigma^{*}_{ij}$.
Let
$$
  \Omega^{'} = \{ s \in \L_2  ~|~ | \d_{\L_1+s}( E_2) - \beta_2 | \ge 4\a_0^{1/2} \mbox{ or }
$$
$$
                                    \frac{1}{|\L_1|} \sum_{n\in \L_1 +s} |\d_{\L'+n} (E_2) - \beta_2 |^2  \ge 4 \a_0^{1/2}  \},
  \mbox{ and }
  G^{'} = \L_2 \setminus \Omega^{'} \,.
$$
By assumption $\L_1$ is $\eps_0$ attendant of $\L_2$
and $E_2$ is $(\a_0,\eps)$--uniform subset of $\L_2$.
Using Lemma \ref{l:intermediate}, we get $|\Omega^{'}| \le 8 \a_0^{1/2} |\L_2|$.
Let $\tilde{\L} = \L' + \L' + \L' + \L' + \L_1$.
Since  $\L'$ is $\eps$  attendant of $\L_1$, it follows that for any $s\in G^{'}$ we have
$ | \d_{ \tilde{\L} + s } ( E_2) - \beta_2 | < 8 \a_0^{1/2}$
and
$ \sum_{n\in \tilde{\L} +s} |\d_{\L'+n} (E_2) - \beta_2 |^2  < 8 \a_0^{1/2} |\t{\L}|$.
For an arbitrary  $i\in \L_1$ consider the set
$$
 \Omega^{*} = \Omega^{*}_i = \{~ j \in \L_2 ~|~ | \d_{ \tilde{\L}_i + j } ( E_2) - \beta_2 | \ge 8 \a_0^{1/2}
                                                    \mbox{ or }
$$
\begin{equation}\label{f:tilde}
                                    \frac{1}{|\t{\L}_i|} \sum_{n\in \t{\L}_i +j} |\d_{\L'+n} (E_2) - \beta_2 |^2  \ge 8\a_0^{1/2} ~\} \,.
\end{equation}
Since  $ (\L_2 \setminus \Omega^{*}_i ) \supseteq ( \L_2 \cap  (G^{'} - i) )$,
it follows that $\Omega^{*}_i \subseteq ( \L_2 \setminus (G^{'} - i) )$.
Since  $\L_1$ is $\eps_0$ attendant of $\L_2$, it follows that
$| \L_2 \setminus (G^{'} - i) | = | (\L_2 + i) \setminus G^{'} | \ge |\L_2^{-} \cap G^{'}|
\ge (1 - 8 \a_0^{1/2} - 8\k_0) |\L_2|$, $\k_0 \le \a_0^2$.
Hence $|\Omega^{*}_i| \le 8 \a_0^{1/2} |\L_2| + 8 \k_0 |\L_2| \le 16 \a_0^{1/2} |\L_2|$.
This yields
\begin{equation}\label{est:*}
  \frac{1}{|\L'|^2}
    \sum_{i\notin \Omega_1, j \in \Omega^{*}_i} |\sigma_{ij}|
        \le
            \frac{1}{|\L'|^2}  \sum_{i\notin \Omega_1} ( 16 \a_0^{1/2} |\L_2|  2 |\L'|^2 |\L_1| )
        \le
            32 \a_0^{1/2} |\L_1|^2 |\L_2| \,.
\end{equation}

We have $j \notin \Omega_2$.
Suppose in addition that  $j\notin \Omega^{*}_i$.
Using   (\ref{svertka}), we get
\begin{equation}\label{e:tt_1}
  \sigma^{*}_{ij} \le \int_0^1 | \F{E}_1^{(2)} (x) |^2 | \F{\tilde{E}}_2^{(2)} (x) |^2 dx \,.
\end{equation}
Since  $E_1^{(2)}$ is $\a_0$--uniform, it follows that
\begin{equation}\label{e:tt_2}
  \F{E}_1^{(2)} (x) = \beta_1^{(2)} \F{\L}_1^{(2)} (x)  +  \vartheta_1 \a_0 |\L'| \,,
\end{equation}
where $|\vartheta_1| \le 1$.
We have $| \F{\L}_1^{(2)} (x) | \le |\L'|$.
Combining (\ref{e:tt_2}) and (\ref{e:tt_1}), we obtain
$$
  \sigma^{*}_{ij}
                 \le ( \beta_1^{(2)} )^2 \int_0^1 | \F{\L}_1^{(2)} (x) |^2 | \F{\tilde{E}}_2^{(2)} (x) |^2 dx
                 + 3 \a_0 |\L'|^2 \int_0^1 | \F{\tilde{E}}_2^{(2)} (x) |^2 dx \le
$$
\begin{equation}\label{}
             \le ( \beta_1^{(2)} )^2 \sum_m ( \L_1^{(2)} * \tilde{E}_2^{(2)} )^2 (m)
                 + 15 \a_0 |\L'|^2 |\L_1| \,.
\end{equation}

Let
$\t{\beta}_2^{j} = \d_{Q} ( \t{E}_2^{(2)} )$.
Since  $j\notin \Omega^{*}_i$, it follows that $| \t{\beta}_2^{j} - \beta_2 | \le 8 \a_0^{1/2}$
and
$$
  \sum_m ( \L_1^{(2)} * \tilde{E}_2^{(2)} )^2 (m) \le 4 ( \beta_1^{(2)} )^2 \beta_2^2 |\L'|^2 |\L_1| + 200 \a_0^{1/2} |\L'|^2 |\L_1| \,.
$$
This implies that
$$
  \sigma^{*}_{ij} \le 4 ( \beta_1^{(2)} )^2 \beta_2^2 |\L'|^2 |\L_1| + 200 \a_0^{1/2} |\L'|^2 |\L_1| + 15 \a_0 |\L'|^2 |\L_1|
             \le
$$
\begin{equation}\label{}
             \le
             4 ( \beta_1^{(2)} )^2 \beta_2^2 |\L'|^2 |\L_1| + 256 \a_0^{1/2} |\L'|^2 |\L_1| \,.
\end{equation}
Since  $i \notin \Omega_1$, it follows that $\beta_1 /2 \le \beta_1^{(2)}  \le 2 \beta_1 $.
Hence $256 \a_0^{1/2} \le (\beta_1^{(2)} )^2 \beta_2^2$.
Consequently for all $i\notin \Omega_1$, $j \notin \Omega_2 \cup \Omega^{*}_i$, we obtain
\begin{equation}\label{}
  \sigma^{*}_{ij} \le 128 ~ \beta_1^2 \beta_2^2 |\L'|^2 |\L_1| \,.
\end{equation}
\\

We have
\begin{equation}\label{}
  |\sigma|^4 \le \| h \|_2^4 \cdot \sigma^{*} \cdot \sum_{m,u} \sum_{r,r'} f(r,m) f(r,u) f(r',m) f(r',u) \,\, \cdot
\end{equation}
\begin{equation}\label{}
                                      \sum_k E_1^{(2)} (k) E_1^{(2)} (k-m+u)
                                             E_2^{(1)} (m-k+r) E_2^{(1)} (m-k+r') =
\end{equation}
\begin{equation}\label{}
           =   \| h \|_2^4 \cdot \sigma^{*} \cdot \sum_{m,u} \sum_{r,r'} f(r,m) f(r,u) f(r',m) f(r',u) \,\, \cdot
\end{equation}
\begin{equation}\label{e:final_E}
                                      \sum_k E_1^{(2)} (m-k) E_1^{(2)} (u-k)
                                             E_2^{(1)} (k+r) E_2^{(1)} (k+r')
           = \| h \|_2^4 \cdot \sigma^{*} \cdot \sigma' \,.
\end{equation}
Rewrite  $\sigma'$ as
\begin{equation}\label{e:final_E1}
  \sigma' = \sum_k \sum_{r,r'} E_2^{(1)} (k+r) E_2^{(1)} (k+r')
                                                  | \sum_m E_1^{(2)} (m-k) f(r,m) f(r',m) |^2
\end{equation}
We have  $r\in \L_1$ and $k+r \in \L_2^{(1)}$. It follows that $k\in \L_2^{(1)} - \L_1$.
On the other hand $m-k \in \L_1^{(2)}$.
Hence $m\in \L_1^{(2)} + k \in \L_2^{(1)} + \L_1^{(2)} - \L_1$.
By symmetry  $u$ belongs to $\L_2^{(1)} + \L_1^{(2)} - \L_1$.
Using Lemma \ref{l:L_pm} for $\L_1$ and its $\eps$ attendant
$\L'$, we obtain that $k$ and $m,u$ belongs to
some translations of {\it Bohr } sets
$W_1 = \L_1^{+}$ and $W_2 = W_1^{+}$, respectively,
and
the cardinalities of these sets do not exceed $3|\L_1|$.

If $k$ is fixed, then $m,u,r,r'$ in (\ref{e:final_E})
run some sets of the cardinalities at most $|\L'|$.

Let $\Phi^1_{r,r'} (m) = f(r,-m) f(r',-m) W_2(m-i-j)$,\\
$\Phi^2_{r,r'} (u) = f(r,-u) f(r',-u) W_2 (u-i-j)$,
$\Phi^3_{m,u} (r) = f(-r,m) f(-r,u)$, and $\Phi^4_{m,u} (r') = f(r',m) f(r',u)$.
Consider the sets
$$
 B_1 = \{ k ~|~ | (\Phi^1_{r,r'} * E_1^{(2)}) (-k) - \beta_1^{(2)} (\Phi^1_{r,r'} * \L_1^{(2)}) (-k) | \ge \a_0^{2/3} |\L'| \} \,
$$
$$
 B_2 = \{ k ~|~ | (\Phi^2_{r,r'} * E_1^{(2)}) (-k) - \beta_1^{(2)} (\Phi^2_{r,r'} * \L_1^{(2)}) (-k) | \ge \a_0^{2/3} |\L'| \} \,
$$
$$
 B_3 = \{ k \in \L_1 ~|~ | (\Phi^3_{m,u} * E_2^{(1)}) (k) - \beta_2^{(1)} (\Phi^3_{m,u} * \L_2^{(1)}) (k) | \ge \a_0^{2/3} |\L'| \} \,
$$
$$
 B_4 = \{ k \in \L_1 ~|~ | (\Phi^4_{m,u} * E_2^{(1)}) (k) - \beta_2^{(1)} (\Phi^4_{m,u} * \L_2^{(1)}) (k) | \ge \a_0^{2/3} |\L'| \} \,.
$$
We have $i \notin \Omega_1$, $j \notin \Omega_2$.
Using Corollary  \ref{c:sv}, we get
$|B_1|, |B_2| \le 3 \a_0^{2/3} |\L_1|$
and
$|B_3|, |B_4| \le \a_0^{2/3} |\L_1|$.
Let $B = B_1 \bigcup B_2 \bigcup B_3 \bigcup B_4$.
Then $|B| \le 8 \a_0^{2/3} |\L_1|$.
Split $\sigma'$ as
$$
  \sigma' = \sum_{k\in B} \sum_{r,r'} E_2^{(1)} (k+r) E_2^{(1)} (k+r') | \sum_m E_1^{(2)} (m-k) f(r,m) f(r',m) |^2 +
$$
$$
           + \sum_{k \notin B} \sum_{r,r'} E_2^{(1)} (k+r) E_2^{(1)} (k+r') | \sum_m E_1^{(2)} (m-k) f(r,m) f(r',m) |^2
           = \sigma_1 + \sigma_2
$$
Let us estimate $\sigma_1$.
Since  $|B| \le 8 \a_0^{2/3} |\L_1|$, it follows that
\begin{equation}\label{}
  |\sigma_1| \le 8 \a_0^{2/3} |\L'|^4 |\L_1| \,.
\end{equation}
If $k\notin B$, then $k\notin B_1$.
This implies that
$$
  \sigma_2 = \sum_{k\notin B} \sum_{u} \sum_{r,r'} f(r,u) f(r',u) E_1^{(2)} (u-k) E_2^{(1)} (k+r) E_2^{(1)} (k+r') \,\, \cdot
$$
$$
             \sum_m f(r,m) f(r',m) E_1^{(2)} (m-k) =
$$
$$
           = \sum_{k\notin B} \sum_{u} \sum_{r,r'} f(r,u) f(r',u) E_1^{(2)} (u-k) E_2^{(1)} (k+r) E_2^{(1)} (k+r')
            (\Phi^1_{r,r'} * E_1^{(2)}) (-k)
$$
$$
           = \beta_1^{(2)} \sum_{k\notin B} \sum_{u} \sum_{r,r'} f(r,u) f(r',u) E_1^{(2)} (u-k) E_2^{(1)} (k+r) E_2^{(1)} (k+r')
             \,\, \cdot
$$
$$
                                                \sum_m f(r,m) f(r',m) \L_1^{(2)} (m-k) +
$$
$$
           + \vartheta \a_0^{2/3} |\L'| \sum_{k\notin B} \sum_{u} \sum_{r,r'} f(r,u) f(r',u) E_1^{(2)} (u-k) E_2^{(1)} (k+r) E_2^{(1)} (k+r')
$$
$$
           = \beta_1^{(2)} \sum_{k\notin B} \sum_{u} \sum_{r,r'} f(r,u) f(r',u) E_1^{(2)} (u-k) E_2^{(1)} (k+r) E_2^{(1)} (k+r') \,\, \cdot
$$
\begin{equation}\label{}
                                                \sum_m f(r,m) f(r',m) \L_1^{(2)} (m-k)
           + 4 \vartheta \a_0^{2/3} |\L'|^4 |\L_1| \,,
\end{equation}
where $|\vartheta| \le 1$.
Using these arguments for $B_2$, $B_3$ and $B_4$, we get
$$
  |\sigma_2| \le (\beta_1^{(2)})^2 (\beta_2^{(1)})^2
                                      \sum_{m,u} \sum_{r,r'} f(r,m) f(r,u) f(r',m) f(r',u) \,\, \cdot
$$
\begin{equation}\label{}
                                      \sum_k \L_1^{(2)} (m-k) \L_1^{(2)} (u-k)
                                             \L_2^{(1)} (k+r) \L_2^{(1)} (k+r')
                                      +
                                      16 \a_0^{2/3} |\L'|^4 |\L_1| \,,
\end{equation}
It follows that
$$
  |\sigma'| \le |\sigma_1| + |\sigma_2|
                                      \le (\beta_1^{(2)})^2 (\beta_2^{(1)})^2
                                      \sum_{m,u} \sum_{r,r'} f(r,m) f(r,u) f(r',m) f(r',u) \,\, \cdot
$$
\begin{equation}\label{}
                                      \sum_k \L_1^{(2)} (m-k) \L_1^{(2)} (u-k)
                                             \L_2^{(1)} (k+r) \L_2^{(1)} (k+r')
                                      +
                                      32 \a_0^{2/3} |\L'|^4 |\L_1| \,.
\end{equation}
Using (\ref{e:final_E}), we obtain
$$
  |\sigma|^4 \le \| h \|_2^4 \cdot \sigma^{*} \cdot (\beta_1^{(2)})^2 (\beta_2^{(1)})^2
                                                  \sum_k \sum_{r,r'} \L_2^{(1)} (k+r) \L_2^{(1)} (k+r') \,\, \cdot
$$
\begin{equation}\label{e:final_E_f}
                                                  \Big| \sum_m \L_1^{(2)} (m-k) f(r,m) f(r',m) \Big|^2
                                                  +
                                                  32 \| h \|_2^4 \cdot \sigma^{*} \cdot \a_0^{2/3} |\L'|^4 |\L_1|
\end{equation}

We have
\begin{equation}\label{}
  \| h \|_2^2 = \sum_{k,m} h(k,m) \le \sum_{k,m} E_1^{(1)} (k) E_2^{(1)} (m) = \beta_1 \beta_2^{(1)} |\L'| |\L_1| \,.
\end{equation}
Since  $i \notin \Omega_1$, $j \notin \Omega_2$, it follows that $\beta_1^{(2)} \le 2\beta_1$ and $\beta_2^{(1)} \le 2 \beta_2$.
Combining the estimates of $\| h \|_2^2$ and $\sigma^{*}$  with (\ref{e:final_E_f}), we get
$$
  |\sigma_{ij}|^4 \le 2^{15} \beta_1^6 \beta_2^6 |\L'|^4 |\L_1|^3
  \sum_k \sum_{r,r'} \L_2^{(1)} (k+r) \L_2^{(1)} (k+r') \, \cdot \,
$$
\begin{equation}\label{}
  \Big| \sum_m \L_1^{(2)} (m-k) f(r,m) f(r',m) \Big|^2 + 2^{15} \a_0^{2/3} |\L'|^8 |\L_1|^4 \,.
\end{equation}
Let $\Omega'_2 = \Omega'_2 (i) = \Omega_2 \cup \Omega^{*}_i$.
We have
$$
  \Big( \sum_{i\notin \Omega_1, j \notin \Omega'_2 (i)} |\sigma_{i,j}| \Big)^4
        \le
           ( |\L_1| |\L_2| )^3 \sum_{i\notin \Omega_1, j \notin \Omega'_2 (i)} |\sigma_{i,j}|^4
        \le
$$
$$
        \le
           2^{15} \beta_1^6 \beta_2^6 |\L'|^4 |\L_1|^3  (|\L_1| |\L_2|)^3
           \sum_{i\notin \Omega_1, j \notin \Omega'_2 (i)}
           \sum_k \sum_{r,r'} \mu_j (k+r) \mu_j (k+r') \, \cdot \,
$$
$$
           \Big| \sum_m \l_i (m-k) f(r,m) f(r',m) \Big|^2
        +
           2^{15} \a_0^{2/3} (|\L_1| |\L_2|)^4 |\L'|^8 |\L_1|^4
        \le
$$
$$
        \le
           2^{15} \beta_1^6 \beta_2^6 |\L'|^4 |\L_1|^3  (|\L_1| |\L_2|)^3
           \sum_{i\in \L_1, j \in \L_2}
           \sum_k \sum_{r,r'} \mu_j (k+r) \mu_j (k+r') \, \cdot \,
$$
$$
           \Big| \sum_m \l_i (m-k) f(r,m) f(r',m) \Big|^2
        +
           2^{15} \a_0^{2/3} (|\L_1| |\L_2|)^4 |\L'|^8 |\L_1|^4 \,.
$$
By assumption the function $f$ is rectilinearly $(\a,\eps)$--uniform.
It follows that
$$
  \Big( \sum_{i\notin \Omega_1, j \notin \Omega'_2 (i) } |\sigma_{i,j}| \Big)^4
        \le
            2^{15} \a \beta_1^8 \beta_2^8 |\L'|^8 |\L_1|^8 |\L_2|^4
        +
            2^{15} \a_0^{2/3} |\L'|^8|\L_1|^8 |\L_2|^4
        \le
$$
\begin{equation}\label{}
        \le
            2^{16} \a \beta_1^8 \beta_2^8 |\L_1|^8 |\L_2|^4 \,.
\end{equation}
Hence
\begin{equation}\label{}
  \sum_{i\notin \Omega_1, j \notin \Omega'_2 (i) } |\sigma_{i,j}| \le 2^4 \a^{1/4} \beta_1^2 \beta_2^2 |\L_1|^2 |\L_2|  \,.
\end{equation}
Using  (\ref{e:sigma'}), (\ref{est:*}) and (\ref{f:4+R}), we have
$$
  |\sigma_0| \le 16 \k |\L_1|^2 |\L_2| + 8 \a_0^{1/2} |\L_1|^2 |\L_2| +
                 32 \a_0^{1/2} |\L_1|^2 |\L_2| + 2^4 \a^{1/4} \beta_1^2 \beta_2^2 |\L_1|^2 |\L_2|
             \le
$$
\begin{equation}\label{}
             \le 2^5 \a^{1/4} \beta_1^2 \beta_2^2 |\L_1|^2 |\L_2|
\end{equation}
as required.

  The next result is the main in this section.

Let $\L_1$, $\L_2$ be Bohr sets, $\L_1 \le \L_2$,
$\L_1 = \L_{\theta,\eps_1, N_1}$,
$\theta \in {\bf T}^d$,
and $E_1 \subseteq \L_1$, $E_2 \subseteq \L_2$,
$|E_1| = \beta_1 |\L_1|$, $|E_2| = \beta_2 |\L_2|$.
Let $\mathcal{P}$ be a product set $E_1 \m E_2$.                                                                

\Th {\it
Let  $A$ be an arbitrary subset of $E_1\m E_2$
of cardinality $\delta |E_1||E_2|$.
Suppose that  the sets $E_1,E_2$ are $(\a_0,2^{-10} \eps^2)$--uniform,
$\a_0 = 2^{-2000} \d^{96} \beta_1^{48} \beta_2^{48}$, $\eps = (2^{-100}  \a_0^2 ) / (100 d)$.
Let  $A$ be rectilinearly $(\a,\a_1,\eps)$--uniform,
$\a = 2^{-100} \d^{12}$, $\a_1 = 2^{-7} \d$, and
\begin{equation}\label{COND:N}
   \log N_1 \ge 2^{10} d \log \frac{1}{\eps_1 \eps} \,.
\end{equation}
Then $A$ contains a triple $\{ (k,m), (k+d,m), (k,m+d) \}$, where $d\neq 0$.
}
\label{a_case}
\\
\Proof
Let $\L'$ be $\eps$ attendant set of $\L_1$, and $\l_i = \L' + i$, $i \in \L_1$.
Let $G_i = (\l_i \m \L_2) \cap A$,
$f_i (\v{s}) = f(s_1 + i, s_2) \L'(s_1, s_2) $, $i\in \L_1$.
By $G_i$ denote the characteristic functions of the sets $G_i$.
Let
$$
  B_1 = \{ i \in \L_1 ~|~ E_1 \cap \l_i  \mbox { is not } (8 \a_0^{1/4}, \eps) \mbox{--uniform}  \} \,,
$$
$$
  B_2 = \{ i \in \L_1 ~|~  | \d_{\l_i} (E_1) - \beta_1 | \ge 4 \a_0^{1/2}  \} \,,
$$
$$
  B_3 = \{ i \in \L_1 ~|~ \| f_i \|^4_{\L' \m \L_2,\eps} >
         \a \beta_1^2 \beta_2^2 |\L'(\eps)|^4 |\L'|^2 |\L_2| \},
  \mbox{ and } B = B_1 \cup B_2 \cup B_3 \,.
$$
By assumption $E_1$ is $(\a_0,\eps)$--uniform.
By Lemma  \ref{l:intermediate}, we get
$|B_1| \le 8 \a_0^{1/4} |\L_1|$.
and $|B_2| \le 8 \a_0^{1/4} |\L_1|$.
Since  $A$ is rectilinearly $(\a,\a_1, \eps)$--uniform,  it follows that $|B_3| \le \a_1 |\L_1|$.
Hence $|B| \le 16 \a_0^{1/4} |\L_1| + \a_1 |\L_1| \le 2 \a_1 |\L_1|$.

Let $\v{e} = \v{e}_1 + \v{e}_2$, and $\v{s} = x \v{e}_1 + y \v{e}$.
Using Lemma \ref{l:L_pm}, we obtain
\begin{equation}\label{tmp:20:08:1}
  A(\v{s}) = \frac{1}{|\L'|} \cdot \sum_{i\in \L_1} G_i (\v{s}) + \epsilon (\v{s}) \,,
\end{equation}
where
$\k =  \a_0^2$.
Сonsider the sum
\begin{equation}\label{}
  \sigma = \frac{1}{|\L'|} \sum_{i\in \L_1} \sum_{\v{s}} G_i (\v{s}) \,.
\end{equation}
We have $|A| = \d \beta_1 \beta_2 |\L_1| |\L_2|$.
Using  (\ref{tmp:20:08:1}), we get
\begin{equation}\label{tmp_14:36.0}
\sigma \ge \frac{\d \beta_1 \beta_2 }{4} |\L_1| |\L_2| \,.
\end{equation}
Split $\sigma$ as
\begin{equation}\label{}
  \sigma = \frac{1}{|\L'|}  \sum_{i\in B} \sum_{\v{s}} G_i (\v{s})
            +
           \frac{1}{|\L'|} \sum_{i\notin B} \sum_{\v{s}} G_i (\v{s})
         = \sigma_1 + \sigma_2 \,.
\end{equation}
Let us estimate $\sigma_1$.
For any $i \in \L_1$ we have $\sum_{\v{s}} G_i (\v{s}) \le \beta_2 |\L'| |\L_2|$.
If $i\notin B_2$, then
$\sum_{\v{s}} G_i (\v{s}) \le 2 \beta_1 \beta_2 |\L'| |\L_2|$.
It follows that
$$
  \sigma_1 \le
                \frac{1}{|\L'|}
                \sum_{i\in B \cap B_2 } \sum_{\v{s}} G_i (\v{s})
            +
                \frac{1}{|\L'|}
                \sum_{i\in B, i \notin B_2 } \sum_{\v{s}} G_i (\v{s})
            \le
$$
\begin{equation}\label{tmp_14:36}
            \le
                8 \a_0^{1/4} |\L_1| |\L_2| + 4 \a_1 \beta_1 \beta_2 |\L_1| |\L_2| \,.
\end{equation}
By assumption $\a_0^{1/4} \le \a_1 \beta_1 \beta_2$.
This implies that
\begin{equation}\label{tmp_20:32}
  \sigma_1 \le 8 \a_1 \beta_1 \beta_2 |\L_1| |\L_2| + 4 \a_1 \beta_1 \beta_2 |\L_1| |\L_2|
                <
                    16 \a_1 \beta_1 \beta_2 |\L_1| |\L_2| \,.
\end{equation}
We have $\a_1 = 2^{-7} \d$.
Using this and (\ref{tmp_14:36.0}), (\ref{tmp_20:32}), we obtain
\begin{equation}\label{tmp_14:40}
  \frac{1}{|\L'|} \sum_{i\notin B} \sum_{\v{s}} G_i (\v{s}) \ge \frac{\d \beta_1 \beta_2 }{8} |\L_1| |\L_2| \,.
\end{equation}
The formula (\ref{tmp_14:40}) implies that there exists $i_0\notin B$ such that
\begin{equation}\label{}
  \sum_{\v{s}} G_{i_0} (\v{s})  \ge \frac{\d \beta_1 \beta_2 }{8} |\L'| |\L_2|
                            = 2^{-3} \d \beta_1 \beta_2 |\L'| |\L_2| \,.
\end{equation}
Let $G(\v{s}) = G_{i_0} (\v{s})$.
We have
\begin{equation}\label{tmp_14:47}
  \sum_k \sum_m G(k,m) \ge 2^{-3} \d \beta_1 \beta_2 |\L'| |\L_2| \,.
\end{equation}
We have $m\in \L_2$ and $k+m \in \l_i$.
It follows that $k\in \l_i - \L_2$.
Using Lemma \ref{l:L_pm} we obtain that $k$ belongs to  a set of cardinality at most $2|\L_2|$.
By the Cauchy--Bounyakovskiy inequality, we get
\begin{equation}\label{}
  2^{-6} \d^2 \beta_1^2 \beta_2^2 |\L'|^2 |\L_2|^2 \le \sum_k \Big( \sum_m G(k,m) \Big)^2 \cdot 2|\L_2| \,.
\end{equation}
It follows that
\begin{equation}\label{e:main_term}
  \sum_k \Big( \sum_m G(k,m) \Big)^2 = \sum_k \sum_{m,p} G(k,m) G(k,p)
                                  \ge 2^{-7} \d^2 \beta_1^2 \beta_2^2 |\L'|^2 |\L_2| \,.
\end{equation}
Let  $\v{s} = x \v{e}_1 + y \v{e}_2$.
Consider the sum
\begin{equation}\label{e:num_corners}
  \sigma_0 = \sum_{\v{s}} \sum_{r} G(\v{s}) G(\v{s} + r \v{e}) A(\v{s} + r \v{e}_2) \,.
\end{equation}
Then
\begin{equation}\label{}
 \sigma_0 = \d \sum_{\v{s}} \sum_{r} G(\v{s}) G(\v{s} + r \v{e}) \mathcal{P} (\v{s} + r \v{e}_2)
            +
            \sum_{\v{s}} \sum_{r} G(\v{s}) G(\v{s} + r \v{e}) f(\v{s} + r \v{e}_2)
\end{equation}
\begin{equation}\label{tmp_15:15}
          = \d \sum_{\v{s}} \sum_{r} G(\v{s}) G(\v{s} + r \v{e})
            +
            \sum_{\v{s}} \sum_{r} G(\v{s}) G(\v{s} + r \v{e}) f(\v{s} + r \v{e}_2)
\end{equation}
Let us estimate the second term in (\ref{tmp_15:15}).
Let $\overline{f} (\v{s}) = f(\v{s})$ if $\v{s} \in \l_i \m \L_2$ and
$\overline{f} (\v{s}) = 0 $ otherwise.
We have
\begin{equation}\label{tmp:kinder_g}
  G(\v{s}) G(\v{s} + r \v{e}) f(\v{s} + r \v{e}_2) = G(\v{s}) G(\v{s} + r \v{e}) \o{f} (\v{s} + r \v{e}_2) \,.
\end{equation}
It follows that
$$
  \sum_{\v{s}} \sum_{r} G(\v{s}) G(\v{s} + r \v{e}) f(\v{s} + r \v{e}_2)
   =
  \sum_{\v{s}} \sum_{r} G(\v{s}) G(\v{s} + r \v{e}) \o{f}(\v{s} + r \v{e}_2)
  =
$$
$$
  =
  \sum_{\v{s}} \sum_{r} G(\v{s} + i_0 \v{e}_1 ) G(\v{s} + r \v{e} + i_0 \v{e}_1 ) f(\v{s} + r \v{e}_2 + i_0 \v{e}_1 )
  =
$$
\begin{equation}\label{}
  =
  \sum_{\v{s}} \sum_{r} G(\v{s} + i_0 \v{e}_1 ) G(\v{s} + r \v{e} + i_0 \v{e}_1 ) f_{i_0}(\v{s} + r \v{e}_2)
\end{equation}
Since  $i_0 \notin B$, it follows that $\| f_{i_0} \|^4 \le \a \beta_1^2 \beta_2^2 |\L'|^2 |\L_2|$
and $\d_{\l_{i_0}} (E_1) \le 2 \beta_1 $.
By assumption $\a = 2^{-100} \d^{12}$.
By Theorem \ref{t:tmp_tha} the second term in (\ref{tmp_15:15})
does not exceed
$
    2^{10} \a^{1/4} \beta_1^2 \beta_2^2 |\L'|^2 |\L_2|
  \le
    2^{-8} \d^3 \beta_1^2 \beta_2^2 |\L'|^2 |\L_2|
$.
The inequality (\ref{e:main_term}) implies that the first term in
(\ref{tmp_15:15}) is greater then
$2^{-7} \d^3 \beta_1^2 \beta_2^2 |\L'|^2 |\L_2|$.
It follows that
$\sigma_0 \ge 2^{-8} \d^3 \beta_1^2 \beta_2^2 |\L'|^2 |\L_2|$.

The sum (\ref{e:num_corners}) is the number of triples
$\{ (k,m), (k+d,m), (k,m+d) \}$, where $k \in \L_{i_0}$, $m\in \L_2$, $d\in {\bf Z}$.
The number of triples with $d=0$ does not exceed $|\L'| |\L_2|$.
By assumption $\log N_1 \ge 2^{10} d \log \frac{1}{\eps_1 \eps} $.
Using Lemma \ref{l:Bohr_est}, we get $|\L'| > 2^8 ( \d^3 \beta_1^2 \beta_2^2 )^{-1}$.
Hence,
$2^{-8} \d^3 \beta_1^2 \beta_2^2 |\L'|^2 |\L_2| > |\L'| |\L_2|$.
It follows that $A$ contains a triple $\{ (k,m), (k+d,m), (k,m+d) \}$ with $d\neq 0$.
This completes the proof.

\refstepcounter{section}

\begin{center}
\textsc{ \arabic{section}. Non--uniform case.}
\end{center}

\Lemma {\it
Let $\L_1$, $\L_2$ be Bohr sets, $\L_1 \le \L_2$, and $\L'$ be $\eps$ attendant set of $\L_1$, $\eps = \k / (100 d)$.
Let set $A$ be a subset of  $C \subseteq \L_1 \m \L_2$
of cardinality $\d |C|$.
By $B$ define the set of $s \in \L_1$ such that
$|A\cap ((\L' + s) \m \L_2) | < (\d - \eta) |C\cap ((\L' + s) \m \L_2) |$, where $\eta > 0$.
Then
$$
  \sum_{s \in (\L_1 \setminus B) } |A\cap ((\L' + s) \m \L_2)| \ge
  \d \sum_{s \in (\L_1 \setminus B)} |C\cap ((\L' + s) \m \L_2)|
        +
$$
$$
        +
  \eta \sum_{s \in B} |C\cap ((\L' + s) \m \L_2)| - 4 \k |\L'| |\L_1| |\L_2| \,.
$$
}
\label{l:easy_case}
\Proof
Let $\v{s} = k \v{e}_1 + m \v{e}_2$.
Using Lemma \ref{l:L_pm}, we get
\begin{equation}\label{tmp_13:12;}
  \d |C| = \sum_{\v{s}} A(\v{s}) \L_1 (k) \L_2(m)  =
            \frac{1}{|\L'|} \sum_{n\in \L_1} \sum_{\v{s}} A(\v{s}) ((\L' + n) \m \L_2) (\v{s}) + 2 \vartheta \k |\L_1| |\L_2| \,,
\end{equation}
where $|\vartheta| \le 1$.
Split the sum (\ref{tmp_13:12;}) into a sum over $n \in B$ and a sum over $n \in \L_1\setminus B$.                
We have
$$
  \d |C| < \frac{1}{|\L'|} (\d - \eta) \sum_{n \in B} |C\cap ((\L' + n) \m \L_2)|
                +
$$
\begin{equation}\label{tmp:11113}
                +
            \frac{1}{|\L'|} \sum_{n \in (\L_1 \setminus B)} |A \cap ((\L' + n) \m \L_2)|
            +
            2 \k |\L_1| |\L_2| \,.
\end{equation}
In the same way
\begin{equation}\label{tmp:11132}
  |C| = \frac{1}{|\L'|} \sum_{n \in B} |C\cap ((\L' + n) \m \L_2)|
                +
            \frac{1}{|\L'|} \sum_{n \in (\L_1 \setminus B)} |C\cap ((\L' + n) \m \L_2)| + 2\vartheta_1 \k |\L_1| |\L_2| \,,
\end{equation}
where $|\vartheta_1| \le 1$.
Combining  (\ref{tmp:11113}) and (\ref{tmp:11132}), we obtain the required result.

 \Pred (B. Green)
{\it Let $A$ be a subset of  $E_1\times E_2$
of cardinality $|A| = \delta |E_1| |E_2|$.
Suppose that $\a >0$ is a real number, and $A$ is not rectilinearly $\a$--uniform.
Then
there are two sets $F_1 \subseteq E_1$ and $F_2 \subseteq E_2$ such that
\begin{equation}\label{conj21}
  |A\bigcap (F_1 \m F_2)| > (\d + 2^{-14} \a^2) |F_1||F_2| \quad
\mbox{ and }
\end{equation}
\begin{equation}\label{conj2}
|F_1| \ge 2^{-8} \a |E_1|\,, \quad |F_2| \ge 2^{-8} \a |E_2| \,.
\end{equation}
}
\label{na_case_pr}


In \cite{Shkr_tri} the author used spectral methods to prove Proposition \ref{na_case_pr}.
His proof gives worse constants than  (\ref{conj21}), (\ref{conj2}).
B. Green \cite{Green_BCC} took a more simple approach, which provided better bounds.

Let $\L_1$, $\L_2$ be Bohr sets,
$\L_1 \le \L_2$,
$\L_1 = \L_{\theta,\eps_0, N}$,
$\theta \in {\bf T}^d$,
and $E_1 \subseteq \L_1$, $E_2 \subseteq \L_2$,
$|E_1| = \beta_1 |\L_1|$, $|E_2| = \beta_2 |E_2|$.
Let $\mathcal{P}$ be a product set $E_1 \m E_2$.

\Th
{\it
Let $A$ be a subset of $\mathcal{P}$
of cardinality $|A|=\delta |E_1||E_2|$.
Suppose that $A$ has no triples $\{ (k,m), (k+d,m), (k,m+d) \}$ with $d\neq 0$,
$E_1,E_2$ are $(\a_0,2^{-10} \eps^2)$--uniform,
$\a_0 = 2^{-2000} \d^{96} \beta_1^{48} \beta_2^{48}$, $\eps = (2^{-100} \a_0^2) / (100 d)$,
$\eps' = 2^{-10} \eps^2$,
and
$$
  \log N \ge 2^{10} d \log \frac{1}{\eps_0 \eps} \,.
$$
Then there exists a Bohr set $\t{\L}$, two sets $F_1$, $F_2$
and a vector $\v{y} = (y_1, y_2) \in {\bf Z}^2$,
$F_1 \subseteq E_1 \cap (\t{\L} + y_1) $, $F_2 \subseteq E_2 \cap (\t{\L} + y_2) $ such that
\begin{equation}\label{est::delta}
  \quad |F_1| \ge 2^{-125} \d^{12} \beta_1 |\t{\L}|,  \quad
        |F_2| \ge 2^{-125} \d^{12} \beta_2 |\t{\L}| \,\, \mbox{ and }
\end{equation}
\begin{equation}\label{est::card}
  \quad  \d_{F_1 \m F_2} (A) \ge \d + 2^{-500} \d^{37} \,.
\end{equation}
  Besides that for
  $\t{\L} = \L_{\t{\theta},\t{\eps},\t{N}}$
  we have
  $\t{\theta} = \theta$, $\t{\eps} \ge 2^{-5} \eps' \eps_0 $ and                                      
  $\t{N} \ge 2^{-5} \eps' N$.
}
\label{t:Phase1}
\\
\Proof
Let $\L'$ be $\eps$ attendant of  $\L_1$, and $\L''$ be  $\eps$ attendant of $\L'$.
Suppose that  $A$ is rectilinearly $(\a,\a_1,\eps)$--uniform,
$\a = 2^{-100} \d^{12}$, $\a_1 = 2^{-7} \d$.
Using Theorem \ref{a_case}, we obtain that $A$ contains a triple
$\{ (k,m), (k+d,m), (k,m+d) \}$ with $d\neq 0$.
Hence, the set $A$ is not rectilinearly $(\a,\a_1,\eps)$--uniform.

Let
$$
  B_1 = \{ s \in \L_1 ~|~ | \d_{\L'+s} (E_1) - \beta_1 | \ge 4\a_0^{1/2} \} \,,
$$
$$
  B_2 = \{ s \in \L_1 ~|~ \L' \cap (E_1 -s )\mbox {  is not } (8\a_0^{1/4},\eps) \mbox{--uniform}  \} \,,
$$
and
$$
  B = \{ i \in \L_1 ~|~ \| f_i \|^4_{\L' \m \L_2,\eps}
                            > \a \beta_1^2 \beta_2^2 |\L'(\eps)|^4 |\L'|^2 |\L_2| \} \,.
$$
Since  $A$ is not rectilinearly $(\a,\eps,\eps')$--uniform,
it follows that $|B| > \a_1 |\L_1|$.
By assumption  $E_1$, $E_2$ are $(\a_0,\eps')$--uniform.
Using Lemma \ref{l:intermediate}, we obtain $|B_1| \le 4\a_0^{1/2} |\L_1|$,
$|B_2| \le 8\a_0^{1/2} |\L_1|$.
Let $B_3 = B_1 \cup B_2$. Then $|B_3| \le 12 \a_0^{1/2} |\L_1|$.
Let $B' = B \setminus B_3$.
Since  $32 \a_0^{1/2} < \a_1$, it follows that $|B'| \ge \a_1 |\L_1| /2$.
Note that for all $l\in B'$ we have
\begin{equation}\label{tmp_21:02}
  | \d_{\L'+s} (E_1) - \beta_1 | <  4\a_0^{1/2} \,.
\end{equation}

Let $\eta = 2^{-100} \a^3$.
Let $\l_l = \L' + l$, $l\in \L_1$.
Suppose that for any $l \in B'$ we have
\begin{equation}\label{}
  | A \cap (\l_l \m \L_2) | \le (\d - \eta) |\l_l \cap E_1| |\L_2 \cap E_2| \,.
\end{equation}
Let $B'^{c} = \L_1 \setminus B'$.
Using Lemma \ref{l:easy_case} and  (\ref{tmp_21:02}), we get
$$
  \sum_{l \in B'^{c}} | A\cap (\l_l \m \L_2)| \ge
                                                \d |\L_2 \cap E_2| \sum_{l\in B'^{c}} |\l_l \cap E_1|
                                                +
                                                \eta |\L_2 \cap E_2| \sum_{l\in B'} |\l_l \cap E_1|
                                                - \a_0^2 |\L'| |\L_1| |\L_2|
                                                  \ge
$$
$$
                                                  \ge
                                                \d \beta_2 |\L_2| \sum_{l\in B'^{c}} |\l_l \cap E_1|
                                                +
                                                \eta \frac{\a_1 |\L_1|}{2} \frac{\beta_1 |\L'|}{4} \beta_2 |\L_2|
                                                   \ge
$$
\begin{equation}\label{tmp_21:14}
                                                   \ge
                                                \d \beta_2 |\L_2| \sum_{l\in B'^{c}} |\l_l \cap E_1|
                                                +
                                                2^{-3} \a_1 \eta \beta_1 \beta_2 |\L'| |\L_1| |\L_2|  \,.
\end{equation}
We have
\begin{equation}\label{tmp_21:39}
  \sum_{l \in B_1} | A\cap (\l_l \m \L_2)| \le 4 \a_0^{1/2} |\L_1| |\L'| |\L_2|
  \le 2^{-4} \a_1 \eta \beta_1 \beta_2 |\L'| |\L_1| |\L_2|  \,.
\end{equation}
Combining (\ref{tmp_21:14}) and (\ref{tmp_21:39}), we obtain
\begin{equation}\label{}
  \sum_{l \in (B'^{c} \setminus B_1)} | A\cap (\l_l \m \L_2)|
    \ge
        \d \beta_1 |\L_2| \sum_{l\in B'^{c}} |\l_l \cap E_1|
    +
        2^{-4} \a_1 \eta \beta_1 \beta_2 |\L'| |\L_1| |\L_2|  \,.
\end{equation}
This implies that, there exists a number  $l \in B'^{c} \setminus B_1$ such that
\begin{equation}\label{tmp:eq_20:24}
  |A \cap (\l_l \m \L_2) | > (\d + 2^{-5} \a_1 \eta) |\l_l \cap E_1| |\L_2 \cap E_2| \,.
\end{equation}
Put $\tilde{\L} = \L'$, $y_1 = l_0$ and $F_1 = (\tilde{\L} + l_0) \cap E_1$.
Since  $l_0 \notin B_1$, it follows that $|F_1| \ge \beta_1 |\tilde{\L}| /2$.
The set $E_2$ is $(\a_0, 2^{-10} \eps^2)$--uniform.
This yields  that there exists a number $a$ such that $F_2 = (\tilde{\L} + a) \cap E_2$
has the cardinality at least $\beta_2 |\tilde{\L}_1| /2$
and for $\v{y} = (l_0,a)$ we have
$$
  |A \cap (\t{\L} + \v{y} ) | > (\d + 2^{-6} \a_1 \eta) |F_1| |F_2| \,.
$$
and the theorem is proven.

Let $\v{x} = r \v{e}_1 + m \v{e}_2$, and $f(\v{x})$ be a balanced function of $A$.
There exists $l_0 \in B'$ such that
$$
  | A \cap (\l_{l_0} \m \L_2) | > (\d - \eta) |\l_{l_0} \cap E_1| |\L_2 \cap E_2| \,.
$$
If
\begin{equation}\label{tmp:20:25}
  | A \cap (\l_{l_0} \m \L_2) | \ge (\d + \eta) |\l_{l_0} \cap E_1| |\L_2 \cap E_2| \,,
\end{equation}
then the theorem is proven.

Hence there exists  $l_0 \in B'$ such that
\begin{equation}\label{e:<eta}
  | \sum_{r,m} f(r,m) \l_{l_0} (r) \L_2 (m) | < \eta |\l_{l_0} \cap E_1| |\L_2 \cap E_2| \,.
\end{equation}

Let $\L_0 = \L'+l_0$.
Put $\nu_i = \L'' + i$, $i\in \L_0$ and $\mu_j = \L'' + j$, $j\in \L_2$.
Consider the sum
\begin{equation}\label{tmp_16:19}
  \sigma^* = \sum_{i\in \L_0} \sum_{j\in \L_2} \sum_k \sum_{m} \sum_{r\in \L_0} f(r,m) \nu_i (m-k) \mu_j (k+r) \,.
\end{equation}
Suppose that $i$ and $j$ are fixed in the sum (\ref{tmp_16:19}).
Using Lemma \ref{l:L_pm}, we obtain that
$k$ runs a set of cardinality at most $2|\L_0|$.
Besides that if $i,j,k$ are fixed, then $m$, $r$ run sets of size at most $|\L''|$.
Using Lemma \ref{l:L_pm} once again, we obtain
\begin{equation}\label{tmp_18:23}
   \sigma^* = |\L''|^2 \sum_k \sum_{m} \sum_{r\in \L_0} f(r,m) \L_0(m-k) \L_2(k+r) + \vartheta \a_0^2 |\L''|^2 |\L_0|^2 |\L_2| \,,
\end{equation}
where $|\vartheta| \le 1$.
Let $\L_3 = \L_2 - \L' - l_0$.
Using Lemma \ref{l:L_pm}, we get $|\L_2| \le |\L_3| \le (1+\a_0^2) |\L_2|$.
Note that $k$ belongs to  the set  $\L_3$  in (\ref{tmp_18:23}).                                                
If $k\in \L_2^{-} - l_0$, then $\L_2(k+r) = 1$, for all $r\in \L_0$.
If $k$ is fixed in (\ref{tmp_18:23}), then $r$ and $m$ run
sets of cardinality at most $|\L_0|$.
It follows that
$$
  \frac{\sigma^*} {|\L''|^2}
            = \sum_{k\in (\L_2^{-} - l_0)} \sum_{m} \sum_{r\in \L_0} f(r,m) \L_0(m-k) +
             \sum_{k\in (\L_3 \setminus (\L_2^{-} - l_0)) } \sum_{m} \sum_{r\in \L_0} f(r,m) \L_0(m-k) \L_2(k+r)
$$
$$
           =
             \sum_{k\in (\L_2^{-} - l_0)} \sum_{m} \sum_{r\in \L_0} f(r,m) \L_0(m-k) + \a_0^2 \vartheta_1 |\L_0|^2 |\L_2|
           =
$$
\begin{equation}\label{}
             \sum_k \sum_{m} \sum_{r\in \L_0} f(r,m) \L_0(m-k) + 2 \a_0^2 \vartheta_2 |\L_0|^2 |\L_2|
           =
             |\L_0| \sum_{m} \sum_{r\in \L_0} f(r,m) + 2 \a_0^2 \vartheta_2 |\L_0|^2 |\L_2| \,,
\end{equation}
where $|\vartheta_1|, |\vartheta_2| \le 1$.
Using (\ref{e:<eta}), we get
\begin{equation}\label{tmp_22:47}
  |\sigma^*| < \eta |\L''|^2 |\L_0| |\L_0 \cap E_1| |\L_2 \cap E_2| + 4 \a_0^2 |\L''|^2 |\L_0|^2 |\L_2|
\end{equation}
If $j$ is fixed, then $k$  runs a set $-\L_0 + j + \L''$ in (\ref{tmp_16:19}).
Clearly, the cardinality of this set does not exceed
$(1 + \a_0^2) |\L'|$.
Hence, replacing $4 \a_0^2 |\L''|^2 |\L_0|^2 |\L_2|$ in (\ref{tmp_22:47})                                               
by $8 \a_0^2 |\L''|^2 |\L_0|^2 |\L_2|$,
we can assume that $k$ runs $-\L_0 + j$ in (\ref{tmp_16:19}).

Since  $l\in B'$, it follows that $\beta_1 |\L_0| /2 \le |\L_0 \cap E_1| \le 2 \beta_1 |\L_0|$.
Besides that $16 \a_0^2 < \eta \beta_1 \beta_2$.
This implies that
$$
  | \sum_{i\in \L_0} \sum_{j\in \L_2} \sum_{k\in -\L_0 + j} \sum_{m} \sum_{r\in \L_0} f(r,m) \nu_i (m-k) \mu_j (k+r) |
    <
$$
\begin{equation}\label{e:<eta2}
        <
            2 \eta |\L''|^2 |\L_0| \cdot |\L_0 \cap E_1| \cdot |\L_2 \cap E_2|
        \le 4 \eta \beta_1 \beta_2 |\L''|^2 |\L_0|^2 |\L_2| \,.
\end{equation}
\\

Let
$$
  \Omega = \{ j \in \L_2 ~|~
              \frac{1}{|\L'|} \sum_{k\in \L'+ j} | \d_{\L''+k} (E_2) - \beta_2 |^2 \ge 4 \a_0^{1/2}  \},
              \mbox{ and }
              G = \L_2 \setminus \Omega \,.
$$
Since  $E_2$ is $(\a_0,\eps')$--uniform, it follows that $|\Omega| \le 8 \a_0^{1/2} |\L_2|$.
Let $i\in \L_0$ be fixed.
Let
$$
  \Omega (i)  = \{ j \in \L_2 ~|~
              \frac{1}{|\L'|} \sum_{k\in -\L_0+ j} | \d_{\L''+i+k} (E_2) - \beta_2 |^2 \ge 4 \a_0^{1/2}  \},
              \mbox{ and }
              G (i) = \L_2 \setminus \Omega(i) \,.
$$
Since
$$
    \sum_{k\in -\L_0+ j} | \d_{\L''+i+k} (E_2) - \beta_2 |^2 = \sum_{k\in \L'+ j + (i - l_0) } | \d_{\L''+k} (E_2) - \beta_2 |^2 \,,
$$
it follows that
$\L_2 \cap ( G + l_0 -i) \subseteq G(i)$.
Hence,
$
 |\Omega(i)| \le |\L_2| - | \L_2 \cap ( G + l_0 - i) |
$.
Since  $i$ belongs to  $\L_0$, this implies that a number $a=l_0 - i$ belongs to  $\L'$.
Using Lemma \ref{l:L_pm} for $\L_2$ and its $\eps$ attendant $\L'$, we get
$(G \cap \L_2^{-}) + a \subseteq \L_2$ and
$$
  |\L_2 \cap (G + a)| \ge |\L_2 \cap ( (G\cap \L_2^{-}) + a)| \ge | (G\cap \L_2^{-}) + a| =  | G\cap \L_2^{-} |
  \ge |G| - 8 \a_0^2 |\L_2| \,.
$$
Hence $|\Omega(i)| \le 8 \a_0^{1/2} |\L_2|$.

Since  $l_0 \in B'$, it follows that
 \begin{equation}\label{tmp:21:20'}
  \frac{1}{|\L'|} \sum_{k\in \L'} | \d_{\L''+k} (E_1 - l_0 \cap \L') - \beta_1 |^2 \le 2^6 \a_0^{1/2}
\end{equation}
It is clear that for any $j$ the sum (\ref{tmp:21:20'}) equals
$$
  \frac{1}{|\L'|} \sum_{k \in -\L_0 + j} |\d_{\L'' + j -k} (E_1 \cap \L_0) - \beta_1 |^2 \,.
$$
Indeed
$$
  \sum_{k \in -\L_0 + j} |\d_{\L'' + j -k} (E_1 \cap \L_0) - \beta_1 |^2
  =
  \sum_{k \in \L' + l_0} |\d_{\L'' + k} (E_1 \cap \L' + l_0) - \beta_1 |^2
  =
$$
$$
  =
  \sum_{k \in \L'} |\d_{\L'' + k} (E_1 - l_0 \cap \L') - \beta_1 |^2
$$

Let
$$
  \Omega_{1} (i,j) = \{ k \in -\L_0 + j ~:~ | \d_{\L''+i+k} (E_2) - \beta_2 | \ge 4 \a_0^{1/8} \} \,,
$$
$$
  \Omega_{2} (i,j) = \{ k \in -\L_0 + j ~:~ |\d_{\L'' + j -k} (E_1 \cap \L_0) - \beta_1 | \ge 4 \a_0^{1/8} \},
  \mbox{ and }
$$
$$
  \Omega_{3} (i,j) = \Omega_1 (i,j) \cup \Omega_2 (i,j) \,.
$$
For all $j\notin \Omega(i)$ we have $|\Omega_{1} (i,j)| \le 2 \a_0^{1/4} |\L'|$.
The inequality (\ref{tmp:21:20'}) implies that $|\Omega_2 (i,j)| \le 4 \a_0^{1/4} |\L'|$.
Hence $|\Omega_3 (i,j) | \le 8 \a_0^{1/4} |\L'|$ if $j\notin \Omega(i)$.
\\

Since  $l_0 \in B'$, it follows that
$$
  \sigma = \sum_{i\in \L_0} \sum_{j\in \L_2} \sum_k \sum_{m,u} \nu_i (m-k) \nu_i (u-k)
             \Big|\sum_r \mu_j (k+r) \t{f}_{l_0} (r,m) \t{f}_{l_0} (r,u) \Big|^2
             \ge
$$
\begin{equation}\label{e:nonunif}
             \ge
                \a \beta_1^2 \beta_2^2 |\L''|^4 |\L_0|^2 |\L_2| \,,
\end{equation}
where $\t{f}_{l_0}$ is a restriction of $f$ to $\l_{l_0} \m \L_2$.
If $j$ is fixed, then  $k$ runs $-\L_0 + j + \L''$  in (\ref{e:nonunif}).
Clearly, the cardinality of this set does not exceed $(1 + \a_0^2) |\L'|$.
Hence, replacing  $\a$ by $\a/2$ in (\ref{e:nonunif}), we can assume that
$k$ runs $-\L_0 + j$ in (\ref{e:nonunif}).
Using $|\Omega (i)| \le 8 \a_0^{1/2} |\L_2|$, we get
$$
  \sigma = \sum_{i\in \L_0} \sum_{j\notin \Omega (i)} \sum_k \sum_{m,u} \nu_i (m-k) \nu_i (u-k)
             \Big|\sum_r \mu_j (k+r) \t{f}_{l_0} (r,m) \t{f}_{l_0} (r,u) \Big|^2
             \ge
$$
\begin{equation}\label{tmp_18:05}
             \ge
                \frac{\a}{4} \beta_1^2 \beta_2^2 |\L''|^4  |\L_0|^2 |\L_2| \,.
\end{equation}
\\                                                                                                

Now we can prove the theorem.

Let
$$
  J = \{
         (i,j,k) ~|~ i\in \L_0,\, j\notin \Omega(i),\, k \notin \Omega_3(i,j)
          \mbox{ such that  }
$$
$$
                                  \sum_{m,u} \nu_i (m-k) \nu_i (u-k)
                                  \Big|\sum_r \mu_j (k+r) \t{f}_{l_0} (r,m) \t{f}_{l_0} (r,u) \Big|^2
                                  \ge
                                  \frac{\a}{64} \beta_1^2 \beta_2^2 |\L''|^4
      \} \,.
$$
Using (\ref{tmp_18:05}), we get
$$
  \sum_{i \in \L_0} \sum_{j \notin \Omega(i)} \sum_{k\notin \Omega_3(i,j)} \sum_{m,u} \nu_i (m-k) \nu_i (u-k)
             \Big|\sum_r \mu_j (k+r) \t{f}_{l_0} (r,m) \t{f}_{l_0} (r,u) \Big|^2
             \ge
$$
\begin{equation}\label{tmp_15:05}
             \ge
                \frac{\a}{8} \beta_1^2 \beta_2^2 |\L''|^4  |\L_0|^2 |\L_2| \,.
\end{equation}

It follows that
$$
  \sum_{(i,j,k) \in J} \sum_{m,u} \nu_i (m-k) \nu_i (u-k)
             \Big|\sum_r \mu_j (k+r) \t{f}_{l_0} (r,m) \t{f}_{l_0} (r,u) \Big|^2
             \ge
$$
\begin{equation}\label{tmp_15:06}
             \ge
                \frac{\a}{16} \beta_1^2 \beta_2^2 |\L''|^4  |\L_0|^2 |\L_2| \,.
\end{equation}
Let us estimate the cardinality of $J$.
For any triple $(i,j,k)$ belongs to $J$ we have
$|E_2 \cap (\nu_i + k) | - \beta_2 |\L''| | \le 4\a_0^{1/8} |\L''|$
and
$|(E_1 \cap \L_0) \cap (\mu_j - k) | - \beta_1 |\L''| | \le 4\a_0^{1/8} |\L''|$.
Using (\ref{tmp_15:06}), we get
\begin{equation}\label{}
  32 |J| \cdot |\L''|^4 \beta_1^2 \beta_2^2                                        
    \ge
  \frac{\a}{16} \beta_1^2 \beta_2^2 |\L''|^4  |\L_0|^2 |\L_2| \,.
\end{equation}
This yields that $|J| \ge 2^{-12} \a |\L_0|^2 |\L_2|$.

Let us assume that for all $(i,j,k) \in J$ we have
\begin{equation}\label{tmp_14:59.0}
  \sum_{m} \sum_{r\in \L_0} f(r,m) \nu_i (m-k) \mu_j (k+r) < - 2^{15} \frac{\eta}{\a} \beta_1 \beta_2 |\L''|^2 \,.
\end{equation}
Using (\ref{e:<eta2}), we get
\begin{equation}\label{}
  \sum_{(i,j,k) \in \o{J} } \sum_{m} \sum_{r\in \L_0} f(r,m) \nu_i (m-k) \mu_j (k+r) \ge
                                                    4 \eta \beta_1 \beta_2 |\L''|^2 |\L_0|^2 |\L_2| \,,
\end{equation}
where $\o{J} = \{ (i,j,k) ~:~ (i,j,k) \in (\L_0 \m \L_2 \m (-\L_0 + j) ) \setminus J \}$.
Since  $|\Omega(i)| \le 8 \a_0^{1/2} |\L_2|$, $i\in \L_0$, it follows that
\begin{equation}\label{}
  \sum_{(i,j,k) \in \o{J}, j \notin \Omega(i)}
  \sum_{m} \sum_{r\in \L_0} f(r,m) \nu_i (m-k) \mu_j (k+r) \ge 2 \eta \beta_1 \beta_2 |\L''|^2 |\L_0|^2 |\L_2| \,.
\end{equation}
Hence, there exist $i$ and $j$, $j \notin \Omega (i)$ such that
\begin{equation}\label{}
  \sum_{k \in Q(i,j)} \sum_{m} \sum_{r\in \L_0} f(r,m) \nu_i (m-k) \mu_j (k+r) \ge \frac{\eta}{2} \beta_1 \beta_2  |\L'|^2 |\L_0| \,,
\end{equation}
where $Q(i,j)$ is a  subset of $-\L_0 + j$.
Since  $j\notin \Omega(i)$, it follows that
$|\Omega_3(i,j)| \le 8 \a_0^{1/4} |\L'|$.
Hence
\begin{equation}\label{tmp_19:59}
  \sum_{k\in Q(i,j) \setminus \Omega_3(i,j) } \sum_{m} \sum_{r\in \L_0} f(r,m) \nu_i (m-k) \mu_j (k+r)
    \ge
  \frac{\eta}{4} \beta_1 \beta_2  |\L''|^2 |\L_0| \,.
\end{equation}
This implies that there exists $k\notin \Omega_3(i,j)$ such that
\begin{equation}\label{tmp_14:50}
  \sum_{m} \sum_{r\in \L_0} f(r,m) \nu_i (m-k) \mu_j (k+r)
    \ge
  \frac{\eta}{8} \beta_1 \beta_2  |\L''|^2 \,.
\end{equation}
Put $\tilde{\L} = \L''$,  $\v{y} = (j-k,k+i)$ and $F_1 = (\tilde{\L} + y_1) \cap (E_1 \cap \L_0)$,
$F_2 = (\tilde{\L} + y_2) \cap E_2$.
Since  $k\notin \Omega_3 (i,j)$, it follows that $\beta_1 |\L''| /2 \le |F_1| \le 2 \beta_1 |\L''|$,
$\beta_2 |\L''| /2 \le |F_2| \le 2 \beta_2 |\L''|$.
Using this and (\ref{tmp_14:50}), we get
$$
  |A \cap (F_1 \m F_2)| = |A \cap ( ( (\mu_j - k) \cap \L_0) \m ((\nu_i + k) \cap \L_2 )  ) |
        \ge
$$
$$
        \ge
            \d |(\mu_j - k) \cap E_1 \cap \L_0| |(\nu_i + k) \cap E_2| + \frac{\eta}{8} \beta_1 \beta_2  |\L''|^2
        \ge
$$
$$
        \ge
            ( \d
            +
            \frac{\eta}{32} ) |F_1| |F_2| \,.
$$
Hence, if for all $(i,j,k) \in J$ we have (\ref{tmp_14:59.0}), then the theorem is proven.

Now assume that there exists a triple
$(i,j,k ) \in J$ such that
\begin{equation}\label{}
  \sum_{m} \sum_{r\in \L_0} f(r,m) \nu_i (m-k) \mu_j (k+r) \ge -2^{15} \frac{\eta}{\a} \beta_1 \beta_2 |\L''|^2 \,.
\end{equation}
We can assume that for all  $(i,j,k) \in J$ we have
\begin{equation}\label{est:den_s}
  |\sum_{m} \sum_{r\in \L_0} f(r,m) \nu_i (m-k) \mu_j (k+r)| \le 2^{15} \frac{\eta}{\a} \beta_1 \beta_2 |\L''|^2 \,.
\end{equation}
Indeed, if
$$
  \sum_{m} \sum_{r\in \L_0} f(r,m) \nu_i (m-k) \mu_j (k+r) > 2^{15} \frac{\eta}{\a} \beta_1 \beta_2 |\L''|^2 \,,
$$
then we might apply the same reasoning as above.
For sets
$\tilde{\L}_1 = \L''$, $\tilde{\L}_2 = \L''$, a vector $\v{y} = (j-k,k+i)$ and $F_1 = (\tilde{\L}_1 + y_1) \cap (E_1 \cap \L_0)$,
$F_2 = (\tilde{\L}_2 + y_2) \cap E_2$
we have
$|F_1| \ge \beta_1 |\tilde{\L}_1| /2$,
$|F_2| \ge \beta_2 |\tilde{\L}_2| /2$
 and
$$
|A \cap (F_1 \m F_2) |
\ge
(\d + 2^{6} \frac{\eta}{\a} ) |F_1| |F_2|.
$$
\\

Since  $(i,j,k) \in J$, it follows that
\begin{equation}\label{tmp_22:06}
  \sum_{m,u \in \nu_i +k }
                          \Big|\sum_{r \in \mu_j - k} \t{f}_{l_0} (r,m) \t{f}_{l_0} (r,u) \Big|^2
                               \ge
                          2^{-6} \a \beta_1^2 \beta_2^2 |\L''|^4 \,.
\end{equation}
Note that $m,u$ belong to $\nu_i+k \cap \L_2$ in (\ref{tmp_22:06})
and $r$ belongs to  a set $\mu_j - k \cap \L_0$.
Put $ \mathcal{L}_1 = \mu_j - k \cap \L_0 $, $ \mathcal{L}_2 = \nu_i + k \cap \L_2 $,
$E_1' = E_1 \cap \mathcal{L}_1$ and  $E_2' = E_2 \cap \mathcal{L}_2$.
We can assume that $\t{f}_{l_0}$ is zero outside
$\mathcal{L}_1 \m \mathcal{L}_2$ in (\ref{tmp_22:06}).
Let $A_1 = A \cap (\mathcal{L}_1 \m \mathcal{L}_2)$,
$\delta_1 = \delta_{E_1' \m E_2'} (A)$,
and $f_1$ be a balanced function of $A_1$.
Using (\ref{est:den_s}), we get $|\d_1 - \d| \le 2^{20} \frac{\eta}{\a} $.
We have $k \notin \Omega_3 (i,j)$.
Using this, we obtain
\begin{equation}\label{}
  \| \t{f}_{l_0} - f_1 \|^4 = |E_1'|^2 |E_2'|^2 ( \delta_1 - \delta )^2 \le
  2^{44} \beta_1^2 \beta_2^2 \frac{\eta^2}{\a^2} |\L'|^4 \,.
\end{equation}
We have $\eta \le 2^{-50} \a^2$.
Using this and Lemma \ref{l:norm}, we get
\begin{equation}\label{}
  \sum_{m,u \in \nu_i +k }
                          \Big|\sum_{r \in \mu_j - k} f_1 (r,m) f_1 (r,u) \Big|^2
                               \ge
                          2^{-7} \a \beta_1^2 \beta_2^2 |\L'|^4 \,.
\end{equation}
Since  $k \notin \Omega_3 (i,j)$, it follows that
$2^{-1} \beta_1 |\L'| \le |E_1'| \le 2 \beta_1 |\L'|$,
$2^{-1} \beta_2 |\L'|  \le |E_2'| \le 2 \beta_2 |\L'|$.
Hence
\begin{equation}\label{}
  \sum_{m,u \in \nu_i +k }
                          \Big|\sum_{r \in \mu_j - k} f_1 (r,m) f_1 (r,u) \Big|^2
                               \ge
                          2^{-11} \a |E_1'|^2 |E_2'|^2 \,.
\end{equation}
Using Proposition \ref{na_case_pr}, we obtain sets $F_1 \subseteq E_1' \subseteq \mu_j - k$,
$F_2 \subseteq E_2' \subseteq \nu_i + k$ such that
$$
  |A \cap (F_1 \m F_2) | \ge
  |A_1 \cap (F_1 \m F_2) | \ge (\d_1 + 2^{-36} \a^2 ) |F_1| |F_2| \ge (\d + 2^{-40} \a^2) |F_1| |F_2|
$$
and
$$
  |F_i| \ge 2^{-19} \a |E_i'| \ge 2^{-25} \a \beta_i |\L'|, \quad i=1,2 \,.
$$
Put $\tilde{\L} = \L''$, $\v{y} = (j-k,k+i)$ and $F_1 = (\tilde{\L}_1 + y_1) \cap (E_1 \cap \L_0)$,
$F_2 = (\tilde{\L}_2 + y_2) \cap E_2$.
The sets $\tilde{\L}$ and $F_1$, $F_2$  satisfy
(\ref{est::delta}), (\ref{est::card}).
This concludes the proof.

\refstepcounter{section}

\begin{center}
\textsc{ \arabic{section}. On dense subsets of Borh sets.}
\end{center}

We need  a simple lemma.

\Lemma
{\it
  Let $\L$ be a Bohr set, $\L'$ be $\eps$ attendant of $\L$, $\eps = \k/(100 d)$,
  and $Q$ be a subset of $\L$.
  Let
  $g : 2^{\bf Z} \m {\bf Z^2} \to {\bf D}$
  be the function
  such that
  $g (\L,\v{x}) = \d^2_{\L+\v{x}} (Q)$.
  Then
  \begin{equation}\label{}
   \frac{1}{|\L|^2} \sum_{\v{x} \in \L} g(\L', \v{x}) \ge g( \L,0) - 8\k \,.
  \end{equation}
}
\label{l:l4}
\Proof
Using the Cauchy--Bounyakovskiy inequality and Lemma \ref{l:smooth_d}, we get
$$
  |\L|^2 \sum_{\v{x} \in \L} g(\L', \v{x})
  \ge
    \Big( \sum_{\v{x} \in \L} \d_{\L' + \v{x}} (Q) \Big)^2
  =
    |\L|^2 (\d_{\L} (Q) + 4 \vartheta \k)^2 \,,
$$
where
$ |\vartheta| \le 1$.
This implies that
$$
  \frac{1}{|\L|^2} \sum_{\v{x} \in \L} g(\L', \v{x}) \ge \d_{\L}^2 (Q) - 8\k = g( \L,0) - 8\k
$$
as required.

\Note
Clearly, the one--dimension analog of Lemma \ref{l:l4} takes place.

\Lemma
{\it
  Let $\L$ be a Bohr set, $\L'$ be $\eps$ attendant of $\L$, $\eps = \k/(100 d)$,
  $\a>0$ be a real number,
  and
  $Q$ be a subset of $\L$, $|Q| = \d |\L|$.
  Suppose that
  \begin{equation}\label{e:Big_sq}
    \frac{1}{|\L|^2} \sum_{\v{n} \in \L} | \d_{\L'+\v{n}} (Q) - \d |^2  \ge \a \,.
  \end{equation}
  Then
  \begin{equation}\label{}
    \sum_{\v{n} \in \L} \d^2_{\L'+\v{n}} (Q) \ge \d^2 + \a - 4 \k \,.
  \end{equation}
}
\label{l:l1}
\Proof
Using (\ref{e:Big_sq}), we have
\begin{equation}\label{tmp:15:25}
  \frac{1}{|\L|^2} \sum_{\v{n} \in \L} \d^2_{\L'+\v{n}} (Q) \ge \frac{2 \d}{|\L|^2}
                                                                    \sum_{\v{n} \in \L} \d_{\L'+\v{n}} (Q) - \d^2 + \a \,.
\end{equation}
The first term in (\ref{tmp:15:25}) equals
$$
  \frac{2\d}{|\L'|^2 |\L|^2}  \sum_{\v{s}}  Q(\v{s}) \sum_{\v{n}} \L(\v{n}) \L' (\v{s} - \v{n})
  = \frac{2\d}{|\L'|^2 |\L|^2}  \sum_{\v{s}}  Q(\v{s}) (\L * \L') (\v{s}) \,.
$$
Using Lemma \ref{l:L_pm}, we obtain
\begin{equation}\label{}
   \frac{1}{|\L|^2} \sum_{\v{n} \in \L} \d^2_{\L'+\v{n}} (Q) \ge \frac{2 \d}{|\L|^2} \sum_{\v{s}}  Q(\v{s}) \L (\v{s}) - \d^2 + \a - 4 \k
  \ge \d^2 + \a - 4 \k \,.
\end{equation}
This completes the proof.

\Note
Clearly, the one--dimension analog of Lemma \ref{l:l1} takes place.

\Cor {\it
  Let $\L$ be a  Bohr set, $\a>0$ be a real number,
  and $E_1$, $E_2$  be sets, $|E_1 \cap \L| = \beta_1 |\L|$, $|E_2 \cap \L| = \beta_2 |\L|$.
  Suppose that either $E_1$ or $E_2$ does not satisfy (\ref{c:sm_sq}).
  Let $\L'$ be an arbitrary  $(2^{-10} \a^2 \beta_1^2 \beta_2^2) /  (100 d)$ attendant set  of $\L$.
  Then
  \begin{equation}\label{}
    \frac{1}{|\L|^2} \sum_{\v{n} \in \L} \d_{\L'+\v{n}}^2 (E_1 \m E_2) \ge  \beta_1^2 \beta_2^2  ( 1 + \frac{\a^2}{2} ) \,.
  \end{equation}
}
\label{cor:1case_c}
\Proof
Let $\v{n} = (x,y)$ and $\k = 2^{-10} \a^2 \beta_1^2 \beta_2^2$.
We have
\begin{equation}\label{tmp:22:45''}
  \frac{1}{|\L|^2} \sum_{\v{n} \in \L} \d_{\L'+\v{n}}^2 (E_1 \m E_2)
  =
  \Big( \frac{1}{|\L|} \sum_{x \in \L} \d_{\L' + x}^2 (E_1) \Big)
  \Big( \frac{1}{|\L|} \sum_{y \in \L} \d_{\L' + y}^2 (E_2) \Big)
\end{equation}
Without loss of generality it can be assumed that $E_1$ does not satisfy (\ref{c:sm_sq}).
Using Lemma \ref{l:l1}, we get
\begin{equation}\label{Tmp:22:471}
  \frac{1}{|\L|} \sum_{x \in \L} \d_{\L' + x}^2 (E_1) \ge \beta_1^2 + \a^2 - 4 \k \,.
\end{equation}
Let us estimate the second factor in (\ref{tmp:22:45''}).
Using Lemma \ref{l:l4}, we obtain
\begin{equation}\label{Tmp:22:47'1}
  \frac{1}{|\L|} \sum_{y \in \L} \d_{\L' + y}^2 (E_2) \ge \beta_2^2 - 8 \k \,.
\end{equation}
Combining (\ref{Tmp:22:471}) and (\ref{Tmp:22:47'1}), we have
$$
   \frac{1}{|\L|^2} \sum_{\v{n} \in \L} \d_{\L'+\v{n}}^2 (E_1 \m E_2)
        \ge
   ( \beta_1^2 + \a^2 - 4 \k ) ( \beta_2^2 - 8 \k) \ge \beta_1^2 \beta_2^2  ( 1 + \frac{\a^2}{2} ) \,.
$$
This concludes the proof.

The following lemma was proven by  J. Bourgain in \cite{Bu}.
We give his proof for the sake of completeness.

\Lemma
{\it
  Let $\L = \L_{\theta, \eps, M}$ be a Bohr set, $\a>0$ be a real number,
  and $Q$ be a set, $|Q \cap \L| = \d |\L|$.
  Suppose that
  \begin{equation}\label{e:Big_Fourier}
    \| (Q \cap \L - \d \L) \F{} ~ \|_{\infty}  \ge \a |\L| \,.
  \end{equation}
  Then there exists a Bohr set  $\L' = \L_{\theta',\eps',N'}$ such that
  $\L'$ is  $\eps_1$ attendant of $\L$,
  $\eps_1 = \frac{\k}{100 d}$, $\k\le \a/32$
  and
  \begin{equation}\label{}
    \frac{1}{|\L|} \sum_{n \in \L} | \d_{\L'+n} (Q) - \d |^2 \ge \frac{\a^2}{4} \,,
  \end{equation}
  $\theta' \in {\bf T}^{d+1}$.
}
\label{l:l3}
\\
\Proof
Let $Q_1 = Q \cap \L$.
Using (\ref{e:Big_Fourier}), we obtain
\begin{equation}\label{}
  | \F{Q}_1 (x_0) - \d \F{\L} (x_0) | \ge \a |\L| \,,
\end{equation}
where $x_0 \in {\bf T}$.
We have $\L = \L_{\theta,\eps,M}$.
Put $\theta' = \theta \cup \{ x_0 \} \in {\bf T}^{d+1}$
and
$$
  \L' = \L_{\theta', \frac{\kappa}{100 d} \eps, \frac{\kappa}{100 d} M}   \,.
$$
Using Lemma \ref{l:L_pm}, we get
$$
  \F{Q}_1 (x_0) = \sum_n Q (n) \L(n) e^{2\pi i n x_0} = \frac{1}{|\L'|} \sum_n (\L * \L') (n) Q (n) e^{2\pi i n x_0}
                                                                + 2 \k \vartheta  |\L| \,,
$$
where $|\vartheta| \le 1$.
We have
$$
  \F{Q}_1 (x_0) = \frac{1}{|\L'|} \sum_m \sum_n \L' (n-m) \L(m) Q(n) e^{2\pi i n x_0} +  2 \k \vartheta |\L|
                =
$$
$$
                =
                  \frac{1}{|\L'|} \sum_m \sum_n \L' (n-m) \L(m) Q(n) e^{2\pi i m x_0} +
$$
$$
                +
                  \frac{1}{|\L'|} \sum_m \sum_n \L' (n-m) \L(m) Q(n) [ e^{2\pi i n x_0} - e^{2\pi i m x_0} ] + 2 \k \vartheta  |\L|
                =
$$
$$
                =
                  \sum_{m \in \L} \d_{\L' + m} (Q) e^{2\pi i m x_0} +
                  O(\frac{1}{|\L'|}  \sum_m \sum_n \L' (n-m) \L(m) Q(n) | e^{2\pi i (n-m) x_0} - 1| )
                +
$$
\begin{equation}\label{tmp:19:09}
                +
                  2 \k \vartheta |\L|
                =
                  \sum_{m \in \L} \d_{\L' + m} (Q) e^{2\pi i m x_0} +
                  (14 \k + 2 \k) \vartheta  |\L|  \,.
\end{equation}
Using (\ref{e:Big_Fourier}) and (\ref{tmp:19:09}), we obtain
\begin{equation}\label{}
  \Big| \sum_{m \in \L} \d_{\L' + m} (Q) e^{2\pi i m x_0} - \d \sum_{m\in \L} e^{2\pi i m x_0} \Big| \ge \frac{\a}{2} |\L| \,.
\end{equation}
Hence
\begin{equation}\label{}
  \sum_{m \in \L} | \d_{\L' + m} (Q)  - \d | \ge \frac{\a}{2} |\L| \,.
\end{equation}
Using the Cauchy--Bounyakovskiy inequality, we get
\begin{equation}\label{}
    \frac{1}{|\L|} \sum_{\v{n} \in \L} | \d_{\L'+\v{n}} (Q) - \d |^2 \ge \frac{\a^2}{4} \,.
\end{equation}
This completes the proof.

\Cor
{\it
  Let $\L$ be a Bohr set, $\a>0$ be a real number,
  and $E_1$, $E_2$ be sets, $|E_1 \cap \L| = \beta_1 |\L|$, $|E_2 \cap \L| = \beta_2 |\L|$
  Suppose that either $E_1$ or $E_2$ satisfies (\ref{e:Big_Fourier}).
  Then there exists $(2^{-10} \a^2 \beta_1^2 \beta_2^2) /  (100 d)$ attendant set
  $\L' = \L_{\theta',\eps',N'}$ of Bohr set $\L$ such that
  \begin{equation}\label{}
    \frac{1}{|\L|^2} \sum_{\v{n} \in \L} \d_{\L'+\v{n}}^2 (E_1 \m E_2) \ge  \beta_1^2 \beta_2^2  ( 1 + \frac{\a^2}{8} )
  \end{equation}
  and
  \begin{equation}\label{f:tmp_star}
    \theta' \in {\bf T}^{d+1} \,.
  \end{equation}
}
\label{cor:together}
\Proof
Let $\v{n} = (x,y)$, and $\k = 2^{-10} \a^2 \beta_1^2 \beta_2^2$.
We have
\begin{equation}\label{tmp:22:45'}
  \frac{1}{|\L|^2} \sum_{\v{n} \in \L} \d_{\L'+\v{n}}^2 (E_1 \m E_2)
  =
  \Big( \frac{1}{|\L|} \sum_{x \in \L} \d_{\L' + x}^2 (E_1) \Big)
  \Big( \frac{1}{|\L|} \sum_{y \in \L} \d_{\L' + y}^2 (E_2) \Big)
\end{equation}
We can assume without loss of generality that $E_1$ satisfies  (\ref{e:Big_Fourier}).
Using Lemma \ref{l:l3} and Lemma \ref{l:l1}, we obtain
\begin{equation}\label{Tmp:22:47}
  \frac{1}{|\L|} \sum_{x \in \L} \d_{\L' + x}^2 (E_1) \ge \beta_1^2 + \frac{\a^2}{4} - 4 \k \,.
\end{equation}
Let us estimate the second term in (\ref{tmp:22:45'}).
Using Lemma \ref{l:l4}, we get
\begin{equation}\label{Tmp:22:47'}
  \frac{1}{|\L|} \sum_{y \in \L} \d_{\L' + y}^2 (E_2) \ge \beta_2^2 - 8 \k \,.
\end{equation}
Combining (\ref{Tmp:22:47}) and (\ref{Tmp:22:47'}), we obtain
$$
   \frac{1}{|\L|^2} \sum_{\v{n} \in \L} \d_{\L'+\v{n}}^2 (E_1 \m E_2)
        \ge
   ( \beta_1^2 + \frac{\a^2}{4} - 4 \k ) ( \beta_2^2 - 8 \k) \ge \beta_1^2 \beta_2^2  ( 1 + \frac{\a^2}{8} ) \,.
$$
This concludes  the proof.

We shall say that the vector $\theta'$ from (\ref{f:tmp_star})
is constructed by Corollary  \ref{cor:together}.

Clearly, all lemmas of this section apply to  translations of Bohr sets.
\\

Let  $\mathbf{\L}$ be a union of a family of  Bohr sets 
$\L_0^*,  \L_1^* (\v{x}_0), \dots,  \L_{n}^* (\v{x}_0,\dots, \v{x}_{n-1})$
and
a sequence of some translations of Bohr sets
$\Lambda_0, \L_1(\v{x}_0), \dots, \L_{n} (\v{x}_0,\dots, \v{x}_{n-1})$
such that
$$
  \L_1 (\v{x}_0) \mbox{ and }  \L_1^* (\v{x}_0)  \mbox{ are defined iff } \v{x}_0 \in \L_0              
$$
$$
  \L_2 (\v{x}_0, \v{x}_1) \mbox{ and }  \L_2^* (\v{x}_0, \v{x}_1) \mbox{ are defined iff } \v{x}_1 \in \L_1(\v{x}_0), \v{x}_0 \in \L_0
$$
$$
  \dots
$$
$$
  \L_n (\v{x}_0, \dots, \v{x}_{n-1} ) \mbox{ and }  \L_n^* (\v{x}_0, \dots, \v{x}_{n-1} ) \mbox{ are defined iff }
$$
\begin{equation}\label{c:L}
  \v{x}_{n-1} \in \L_{n-1} (\v{x}_0, \dots, \v{x}_{n-2}),
  \v{x}_{n-2} \in \L_{n-2} (\v{x}_0, \dots, \v{x}_{n-3}), \dots,
  \v{x}_0 \in \L_0 \,.
\end{equation}

Let   $m\ge 0$ be an integer number and $\mathbf{\L}$ be a family of Bohr sets
satisfies  (\ref{c:L}).
Let  $g : 2^{\bf Z} \m {\bf Z^2} \to
{\bf D}$ be a function.
Let us define the  {\it index} of $g$, respect $\mathbf{\L}$, for all $k = 0,\dots, m$ by
$$
  \ind_k (\mathbf{\L}) (g)
    =
  \frac{1}{|\L_0|^2} \sum_{\v{x}_0 \in \L_0} \frac{1}{|\L_1(\v{x}_0)|^2} \sum_{\v{x}_1 \in \L_1 (\v{x}_0) }
    \dots
$$
\begin{equation}\label{i:k}
  \frac{1}{|\L_k(\v{x}_0, \dots, \v{x}_{k-1}) |^2}
  \sum_{\v{y} \in \L_k (\v{x}_0, \dots, \v{x}_{k-1}) }
  g(\L_k^* (\v{x}_0, \dots, \v{x}_{k-1}) ,\v{y}) \,.
\end{equation}

Let $M_k = M_k (\v{x}_0, \dots, \v{x}_{k-1})$  be the family of sets such that
$M_k (\v{x}_0, \dots, \v{x}_{k-1}) \subseteq \L_k (\v{x}_0, \dots, \v{x}_{k-1})$ for all
$(\v{x}_0, \dots, \v{x}_{k-1})$.
For any $k=0,\dots,m$ by  $\ind_k (\mathbf{\L},M) (g)$ define the following expression
$$
  \ind_k (\mathbf{\L},M) (g) =
  \frac{1}{|\L_0|^2} \sum_{\v{x}_0 \in \L_0} \frac{1}{|\L_1(\v{x}_0)|^2} \sum_{\v{x}_1 \in \L_1 (\v{x}_0) }
    \dots
$$
\begin{equation}\label{i:k_M}
  \frac{1}{|\L_k(\v{x}_0, \dots, \v{x}_{k-1}) |^2}
  \sum_{\v{y} \in M_k (\v{x}_0, \dots, \v{x}_{k-1}) }
  g(\L_k^* (\v{x}_0, \dots, \v{x}_{k-1}) ,\v{y}) \,.
\end{equation}
Clearly, we have $| \ind_k (\mathbf{\L}, M) (g) | \le 1$,
for any natural $k\ge 0$, a family $M_k$ and a  function
$g : 2^{\bf Z} \m {\bf Z^2} \to {\bf D}$.

\Lemma
{\it
  Let $Q$ be a subset of $\L_0 \m \L_0$, and $|Q| = \delta |\L_0|^2$.
  Suppose that $\L_k^* (\v{x}_0, \dots, \v{x}_{k-1} )$ is an arbitrary  $\eps$ attendant of
  $\L_k (\v{x}_0, \dots, \v{x}_{k-1} )$, $\eps = \k / (100 d)$.
  Let $ g(M,\v{x}) = \d_{M + \v{x}} (Q)$.
   Then for all $k=0,\dots, n$ we have
  \begin{equation}\label{}
    \Big| \ind_k (\mathbf{\L}) (g) - \d \Big| \le  4 \k (k+1)   \,.
  \end{equation}
}
\label{l:keps}
\Proof
If $k=0$, then Lemma \ref{l:smooth_d} implies the result.
Let $k>0$.
Using Lemma \ref{l:smooth_d} once again, we get
$$
  \ind_k (\mathbf{\L}) (g) \ge
    \ind_{k-1} (\mathbf{\L}) (g) - 4\k \ge
        \dots \ge
            \ind_0 (\mathbf{\L}) (g) - 4 \k k \ge
                \d - 4 \k (k+1) \,.
$$
In the same way we obtain the reverse inequality.
This completes the proof.

   The next result is the main in this section.

\Pred
{\it
Let $\L= \L(\theta, \eps_0, N)$ be a Bohr set, $\theta \in {\bf T}^d$,
and $\v{s} = (s_1,s_2)$ be an integer vector.
Let $\eps, \sigma,\tau, \d \in (0,1)$ be real numbers,
$E_1$, $E_2$ be sets, $E_i = \beta_i |\L|$, $i=1,2$.
Suppose that  ${\bf E} = E_1 \m E_2$ is a subset of  $(\L+s_1) \m (\L+s_2)$,
$A \subset {\bf E}$, $\d_{{\bf E}} (A) = \d + \tau$,
and
$\eps \le \k/(100 d)$, $ \k = 2^{-100} ( \tau \beta_1 \beta_2 )^5 \sigma^3 $.
Let
\begin{equation}\label{c:Nge}
  N \ge (2^{-100} \eps_0 \eps)^{ - 2^{100} ( (\tau \beta_1 \beta_2)^{-5} \sigma^{-3} + d )^2  }  \,,
\end{equation}
and
$
 ~\sigma \le 2^{-100} \tau \beta_1 \beta_2
$.
Then there exists a Bohr set $\L'= \L(\theta',\eps',N')$, $\theta' \in {\bf T}^D$,
$
  D\le 2^{30} (\tau \beta_1 \beta_2)^{-5} \sigma^{-3} + d
$,
$
  \eps' \ge (2^{-10} \eps)^D \eps_0
$,
$
  N'\ge (2^{-10} \eps)^D N
$
and an integer vector $\v{t} = (t_1,t_2)$ such that
if
$E_1' = (E_1 - t_1) \cap \L'$, $E_2' = (E_2 - t_2) \cap \L'$, ${\bf E'} = E_1' \m E_2'$,
then                                                                                                    
\\
$ 1)~ |{\bf E'}| \ge \beta_1 \beta_2 \tau |\L'| /16$; \\
$ 2)~ E_1', E_2' $ are $(\sigma,\eps)$--uniform subsets of $\L'$; \\
$ 3)~ \d_{ {\bf E'} } (A - \v{t}) \ge \d + \tau/16 .$
}
\label{prop:point_3}
\\
\Proof
Let $\beta = \beta_1 \beta_2$, and $\t{E}_1 = E_1 - s_1$, $\t{E}_2 = E_2 -s_2$, $\t{E} = \t{E}_1 \m \t{E}_2$.
If the sets $\t{E}_1$, $\t{E}_2$ are $(\sigma,\eps)$--uniform subsets of $\L$,
then Proposition \ref{prop:point_3} is proven.

Suppose that
$\t{E}_1$, $\t{E}_2$ are not  $(\sigma,\eps)$--uniform subsets of  $\L$.
We shall construct a family of Bohr sets $\mathbf{\L}$
such that $\mathbf{\L}$ satisfies the
conditions
(\ref{c:L}).
The proof of Proposition \ref{prop:point_3} is a sort of an algorithm.
At the first step of our algorithm we put $\L_0 = \L = \L_{\theta,\eps_0,N}$.
If either $\t{E}_1$ or $\t{E}_2$ does not satisfy (\ref{c:sm_total_F}) with $\alpha = \sigma/2$, then
let $\L_0^*$ be  $\eps$ attendant of  $\L_0$ such that $\L_0^*$ is constructed by Corollary  \ref{cor:together}.        
In the other cases let $\L_0^*$ be $\eps$ attendant of  $\L_0$ with the same  $\theta$.
Define
$$
  R_0 = \{ \v{p} = (p_1,p_2) \in \L_0 ~|~ \t{E}_1 - p_1, \t{E}_2 - p_2 \mbox { are } (\sigma,\eps)\mbox{--uniform in } \L_0^*
$$
$$
           \mbox{ or } \d_{\L_0^* + \v{p} } (\t{E}_1 \m \t{E}_2) < \beta \tau /16
        \}
$$
and $\o{R}_0 = (\L_0 \m \L_0) \setminus R_0$.

Let $\tilde{\L}$ be an arbitrary Bohr set, and $\v{n} \in {\bf Z}^2$ be an arbitrary integer vector.
Put  $g(\tilde{\L}, \v{n}) = \d^2_{ \tilde{\L} + \v{n}} ( \t{E})$,
$g_1 (\t{\L}, \v{x}) = \d_{\t{\L}+\v{n}} (A)$,
$ g_2 (\t{\L},\v{n}) = \d_{\t{E} \cap \t{\L} + \v{n} } (A) $
and
$ g_3(\t{\L}, \v{n}) = \d_{\t{\L} + \v{n}} (\t{E})$.
Clearly,  $g(\tilde{\L}, \v{n}) = g^2_3 (\t{\L}, \v{n})$
and
$ g_1 (\t{\L}, \v{x}) \le g_3 (\t{\L}, \v{n})$.
Besides that, we have
$$
  g_1 (\t{\L},\v{n}) = g_2 (\t{\L}, \v{n}) g_3(\t{\L}, \v{n}) \,.
$$

Let $\mathbf{\L}_0 = \{ \L_0 \}$.
If $\ind_0 (\mathbf{\L}_0,\o{R}_0) (g_3) < \tau \beta /4$, then we stop the algorithm at step $0$.

Using Lemma  \ref{l:smooth_d} and the Cauchy--Bounyakovskiy inequality, we get
\begin{equation}\label{}
  \ind_0 (\mathbf{\L}_0 ) (g) \ge
         \Big( \frac{1}{|\L_0|^2} \sum_{\v{y} \in \L_0} \d_{\L_0^* + \v{y}} (\t{E}) \Big)^2
            \ge
            \beta / 2 \,.
\end{equation}

Let after the $k$th step of the algorithm the family of Bohr sets
$\mathbf{\L}_k$ has been constructed, $k\ge 0$.

Let
$$
  \L_{k+1} (\v{x}_0, \dots, \v{x}_{k} ) = \L_k^* (\v{x}_0, \dots, \v{x}_{k-1}) + \v{x}_k \,, \,
   \v{x}_k \in \L_k (\v{x}_0, \dots, \v{x}_{k-1})
$$
Let $\v{x}_k = (a,b)$, and $\L_k^* = \L_k^* (\v{x}_0, \dots, \v{x}_{k-1})$.
If either $(\t{E}_1 - a) \cap \L_k^* $ or $(\t{E}_2 - b) \cap \L_k^*$ does not satisfy  (\ref{c:sm_total_F})
with $\a = \sigma/2$, then let $\L_{k+1}^* (\v{x}_0, \dots, \v{x}_{k} )$ be  $\eps$ attendant of
$\L_{k}^* (\v{x}_0, \dots, \v{x}_{k} )$
such that
$\L_{k+1}^* (\v{x}_0, \dots, \v{x}_{k} )$
is constructed by Corollary  \ref{cor:together}.
In the other cases let $\L_{k+1}^* (\v{x}_0, \dots, \v{x}_{k} )$ be $\eps$ attendant of
$\L_{k}^* (\v{x}_0, \dots, \v{x}_{k} )$
with the same generative vector.

By $R_{k+1} (\v{x}_0, \dots, \v{x}_{k})$,
$\o{R}_{k+1} (\v{x}_0, \dots, \v{x}_{k})$
denote the sets
$$
  R_{k+1} (\v{x}_0, \dots, \v{x}_{k}) =
            \{ \v{p} = (p_1,p_2) \in \L_k^* (\v{x}_0, \dots, \v{x}_{k-1}) + \v{x}_k
            ~|~ \t{E}_1 - p_1, \t{E}_2 - p_2
$$
$$
   \mbox { are }
   (\sigma,\eps)\mbox{--uniform in } \L_{k+1}^* (\v{x}_0, \dots, \v{x}_{k} )
$$
$$
           \mbox{ or } \d_{\L_{k+1}^* (\v{x}_0, \dots, \v{x}_{k} ) + \v{p} } (\t{E}_1 \m \t{E}_2) < \tau \beta /16
        \}
$$
and
$
  \o{R}_{k+1} (\v{x}_0, \dots, \v{x}_{k}) =
        ( \L_k^* (\v{x}_0, \dots, \v{x}_{k-1}) + \v{x}_k ) \setminus R_{k+1} (\v{x}_0, \dots, \v{x}_{k})
$.

By $E_{k} (\v{x}_0, \dots, \v{x}_{k-1})$ denote the sets
$$
  E_{k} (\v{x}_0, \dots, \v{x}_{k-1}) =
$$
$$
  \{
    \v{p} = (p_1,p_2) \in \L_{k-1}^* (\v{x}_0, \dots, \v{x}_{k-2}) + \v{x}_{k-1}
            ~|~
                \d_{\L_{k}^* (\v{x}_0, \dots, \v{x}_{k-1} ) + \v{p} } (\t{E}_1 \m \t{E}_2) < \tau \beta /16
  \}.
$$
Obviously, $E_{k} (\v{x}_0, \dots, \v{x}_{k-1}) \subseteq R_{k} (\v{x}_0, \dots, \v{x}_{k-1})$, $k=0,1,\dots$

Let $\mathbf{\L'}_{k+1} = \{ \L_{k+1} (\v{x}_0, \dots, \v{x}_{k}) \}$, $\v{x}_k \in \L_k (\v{x}_0, \dots, \v{x}_{k-1})$,
and
$\mathbf{\L}_{k+1} = \{ \mathbf{\L}_k , \mathbf{\L'}_{k+1} \}$.
If $\ind_{k+1} (\mathbf{\L}_{k+1}, \o{R}_{k+1}) (g_3) < \tau \beta /4$, then we stop the algorithm at step $k+1$.

Let $\L_{k-1}^* = \L_{k-1}^* (\v{x}_0, \dots, \v{x}_{k-2})$,
and
$ \beta_k' = \d_{\L_{k-1}^*} (\t{E}_1)$,
$ \beta_k'' = \d_{\L_{k-1}^*} (\t{E}_2)$.
Suppose $\v{x}_{k-1} = (a',b')$ belongs to  $\o{R}_{k-1} (\v{x}_0, \dots, \v{x}_{k-2})$.
Note that $\v{x}_{k-1}$ does not belong to $E_{k-1} (\v{x}_0, \dots, \v{x}_{k-2})$.
Let us consider three cases. \\
Case 1 :  either $( \t{E}_1 - a') \cap \L_{k-1}^*$ or
$(\t{E}_2 - b') \cap \L_{k-1}^*$ does not satisfy (\ref{c:sm_B}). \\
Case 2 :  either $( \t{E}_1 - a') \cap \L_{k-1}^*$ or
$(\t{E}_2 - b') \cap \L_{k-1}^*$ does not satisfy (\ref{c:sm_sq}). \\
Case 3 :  either $( \t{E}_1 - a') \cap \L_{k-1}^*$ or
$(\t{E}_2 - b') \cap \L_{k-1}^*$ does not satisfy (\ref{c:sm_total_F}). \\
Note that  $\a$  equals $\sigma$ in all these cases. \\
Let us consider the following situation  :  either $(\t{E}_1 - a') \cap \L_{k-1}^*$ or
$(\t{E}_2 - b') \cap \L_{k-1}^*$ does not satisfy (\ref{c:sm_total_F}) with $\a= 2^{-4} \sigma^{3/2}$.
Let
\begin{equation}\label{}
  S_0 =
       \frac{1}{ |\L_k (\v{x}_0, \dots, \v{x}_{k-1}) |^2 }
       \sum_{\v{y} \in \L_k (\v{x}_0, \dots, \v{x}_{k-1})} g(\L^*_k (\v{x}_0, \dots, \v{x}_{k-1}), \v{y}) \,,
\end{equation}
where
$\L_{k}^* (\v{x}_0, \dots, \v{x}_{k-1} )$
is $\eps$ attendant of  $\L_{k} (\v{x}_0, \dots, \v{x}_{k-1} )$
such that
$\L_{k}^* (\v{x}_0, \dots, \v{x}_{k-1} )$
is constructed by Corollary  \ref{cor:together}.
Using Corollary \ref{cor:together}, we get
$$
   S_0
        \ge
            g ( \L_k (\v{x}_0, \dots, \v{x}_{k-1}), 0) ( 1 + 2^{-11} \sigma^3 )
        =
$$
\begin{equation}\label{}
        =
            g ( \L_{k-1}^* (\v{x}_0, \dots, \v{x}_{k-2}), \v{x}_{k-1} ) ( 1 + 2^{-11} \sigma^3 ) \,.
\end{equation}
Note that in this case, we have $\dim \L^*_k (\v{x}_0, \dots, \v{x}_{k-1}) = \dim \L_k (\v{x}_0, \dots, \v{x}_{k-1}) + 1$.

Suppose that either $(\t{E}_1 - a') \cap \L_{k-1}^*$ or
$(\t{E}_2 - b') \cap \L_{k-1}^*$ does not satisfy (\ref{c:sm_sq}) with $\a = 2^{-4} \sigma^{3/2}$.
Using Corollary \ref{cor:1case_c}, we  obtain
\begin{equation}\label{}
  S_0
        \ge
              g ( \L_{k-1}^* (\v{x}_0, \dots, \v{x}_{k-2}), \v{x}_{k-1} ) ( 1 + 2^{-11} \sigma^3 ) \,.
\end{equation}
In this case, we have $\dim \L^*_k (\v{x}_0, \dots, \v{x}_{k-1}) = \dim \L_k (\v{x}_0, \dots, \v{x}_{k-1})$.

Finally, suppose that either $(\t{E}_1 - a') \cap \L_{k-1}^*$ or $(\t{E}_2 - b') \cap \L_{k-1}^*$
does not satisfy (\ref{c:sm_B}) with $\a = \sigma$.
Note that
$(\t{E}_1 - a') \cap \L_{k-1}^*$
{\it and}
$(\t{E}_2 - b') \cap \L_{k-1}^*$
satisfy (\ref{c:sm_sq}) with $\a = 2^{-4} \sigma^{3/2}$.
Let $\L_{k}^* = \L_{k}^* (\v{x}_0, \dots, \v{x}_{k} )$.
Define
$$
  B_k (\v{x}_0, \dots, \v{x}_{k-1}) = \{ \v{p} = (p_1,p_2) \in \L_k (\v{x}_0, \dots, \v{x}_{k-1})
            ~:~
$$
$$
                \| ( (\t{E}_1 - p_1)  - \beta_k' \L_{k}^* )\F{}~ \|_{\infty} \ge \sigma |\L_{k}^*|
                    \mbox{ or }
                    \| ( (\t{E}_2 - p_2)  - \beta_k'' \L_{k}^* )\F{}~ \|_{\infty} \ge \sigma |\L_{k}^*| \} \,.
$$
We have
\begin{equation}\label{tmp:12:51}
    | B_k (\v{x}_0, \dots, \v{x}_{k-1}) | \ge  \sigma |\L_k (\v{x}_0, \dots, \v{x}_{k-1})|^2 \,.
\end{equation}
Let
$$
  \t{B}_k (\v{x}_0, \dots, \v{x}_{k-1}) = \{ \v{p} = (p_1,p_2) \in B_k (\v{x}_0, \dots, \v{x}_{k-1})
                                                ~:~
$$
$$
                                              | \d_{\L_{k}^*} (\t{E}_1 - p_1) - \beta_k' | \le \sigma / 8
                                                \quad \mbox { and } \quad
                                              | \d_{\L_{k}^*} (\t{E}_2 - p_2) - \beta_k'' | \le \sigma / 8
                                          \} \,.
$$
For all $ \v{p} \in \t{B}_k (\v{x}_0, \dots, \v{x}_{k-1})$, we have
either
$
  (\t{E}_1 - p_1) \cap \L_{k}^*
$
or
$
 (\t{E}_2 - p_2) \cap \L_{k}^*
$
does not $\sigma/2$--uniform.
The sets  $(\t{E}_1 - a') \cap \L_{k-1}^*$ and  $(\t{E}_2 - b') \cap \L_{k-1}^*$
satisfy  (\ref{c:sm_sq}) with $\a$ equals $2^{-4} \sigma^{3/2}$.
This implies that
\begin{equation}\label{tmp:12:51'''}
  | \t{B}_k (\v{x}_0, \dots, \v{x}_{k-1})|
                \ge  \frac{\sigma}{2} |\L_k (\v{x}_0, \dots, \v{x}_{k-1})|^2
\end{equation}

Suppose that
\begin{equation}\label{c:g_3}
  g_3 (\L_{k-1}^*, \v{x}_{k-1}) = \beta_k' \beta_k'' \ge \tau \beta /8 \,.
\end{equation}
It follows from (\ref{c:g_3}) that
\begin{equation}\label{c:g_3'}
  g_3 (\L^*_k (\v{x}_0, \dots, \v{x}_{k-1}), \v{p}) \ge \beta_k' \beta_k'' - \sigma /2 \ge \tau \beta /16 \,,
\end{equation}
for all $\v{p} \in \t{B}_k (\v{x}_0, \dots, \v{x}_{k-1})$.

Let us consider the sum
$$
  S = S (\v{x}_0, \dots, \v{x}_{k-1})
        =
  \frac{1}{ |\L_k (\v{x}_0, \dots, \v{x}_{k-1}) |^2 }
    \sum_{\v{x}_k \in \L_k (\v{x}_0, \dots, \v{x}_{k-1})}
        \frac{1}{ |\L_{k+1} (\v{x}_0, \dots, \v{x}_{k}) |^2 }
$$
$$
    \cdot
    \sum_{\v{y} \in \L_{k+1} (\v{x}_0, \dots, \v{x}_{k})}
    g(\L_{k+1}^* (\v{x}_0, \dots, \v{x}_{k}), \v{y}) \,.
$$
Write the sum $S$ as $S' + S''$, where the summation in $S'$ is taken over
$\v{x}_k \in  \t{B}_k (\v{x}_0, \dots, \v{x}_{k-1})$
and the summation in $S''$ is taken over
$\v{x}_k \in \L_k (\v{x}_0, \dots, \v{x}_{k-1}) \setminus \t{B}_k (\v{x}_0, \dots, \v{x}_{k-1})$.
Note that if $\v{x}_k \in \t{B}_k (\v{x}_0, \dots, \v{x}_{k-1})$, then the Bohr set
$\L_{k+1}^* (\v{x}_0, \dots, \v{x}_{k})$ is constructed by Corollary  \ref{cor:together}.
Using this corollary, we obtain
\begin{equation}\label{tmp:12:52_1}
  S' \ge
    \frac{1}{ |\L_k (\v{x}_0, \dots, \v{x}_{k-1}) |^2 }
        \sum_{\v{y} \in  \t{B}_k (\v{x}_0, \dots, \v{x}_{k-1}) } g(\L_k^* (\v{x}_0, \dots, \v{x}_{k-1}), \v{y}) ( 1 + \frac{\sigma^2}{32})
\end{equation}
Let us estimate the sum $S''$.
Using Lemma \ref{l:l4}, we get
\begin{equation}\label{tmp:12:52_2}
  S'' \ge
        \frac{1}{ |\L_k (\v{x}_0, \dots, \v{x}_{k-1}) |^2 }
            \sum_{\v{y} \in \L_k (\v{x}_0, \dots, \v{x}_{k-1}) \setminus \t{B}_k (\v{x}_0, \dots, \v{x}_{k-1})}
                 g(\L_k^* (\v{x}_0, \dots, \v{x}_{k-1}), \v{y})  - 8 \k
\end{equation}
Combining (\ref{c:g_3'}), (\ref{tmp:12:52_1}), (\ref{tmp:12:52_2})
and (\ref{tmp:12:51'''}), we have
$$
  S \ge
    \frac{1}{ |\L_k (\v{x}_0, \dots, \v{x}_{k-1}) |^2 }
        \sum_{\v{y} \in \L_k (\v{x}_0, \dots, \v{x}_{k-1})} g(\L^*_k (\v{x}_0, \dots, \v{x}_{k-1}), \v{y})
        +
$$
$$
        +
    \frac{1}{ |\L_k (\v{x}_0, \dots, \v{x}_{k-1}) |^2 }
    \sum_{\v{y} \in \t{B}_k (\v{x}_0, \dots, \v{x}_{k-1}) }  2^{-13} \tau^2 \beta^2 \sigma^2   - 2^4 \k
    \ge
$$
$$
    \ge
     \frac{1}{ |\L_k (\v{x}_0, \dots, \v{x}_{k-1}) |^2 }
        \sum_{\v{y} \in \L_k (\v{x}_0, \dots, \v{x}_{k-1})} g(\L^*_k (\v{x}_0, \dots, \v{x}_{k-1}), \v{y})
        +
        2^{-14} \tau^2 \beta^2 \sigma^3 - 2^4 \k
$$
Using Lemma  \ref{l:l4}, we obtain
$$
  S \ge
        g ( \L_{k-1}^* (\v{x}_0, \dots, \v{x}_{k-2}), \v{x}_{k-1} ) + 2^{-14} \tau^2 \beta^2 \sigma^3 - 2^5 \k
    \ge
$$
$$
    \ge
        g ( \L_{k-1}^* (\v{x}_0, \dots, \v{x}_{k-2}), \v{x}_{k-1} ) + 2^{-15} \tau^2 \beta^2 \sigma^3
    \ge
$$
\begin{equation}\label{tmp:12:57}
    \ge
        g ( \L_{k-1}^* (\v{x}_0, \dots, \v{x}_{k-2}), \v{x}_{k-1} ) ( 1 + 2^{-15} \tau^2 \beta^2 \sigma^3 ) \,.
\end{equation}
On the other hand, $S_0$ is an estimate for $S$.
Using Lemma  \ref{l:l4}, we get
$$
  S \ge S_0 - 8 \k \,.
$$
Thus if
$\v{x}_{k-1}$ belongs to  $\o{R}_{k-1} (\v{x}_0, \dots, \v{x}_{k-2})$ and
$\v{x}_{k-1}$ satisfies (\ref{c:g_3}), then we have
\begin{equation}\label{tmp:12:57`}
  S \ge g ( \L_{k-1}^* (\v{x}_0, \dots, \v{x}_{k-2}), \v{x}_{k-1} ) ( 1 + 2^{-15} \tau^2 \beta^2 \sigma^3 ) - 8 \k \,.
\end{equation}
Now suppose that  $\v{x}_{k-1}$ is an arbitrary vector, $\v{x}_{k-1} \in \L_{k-1} (\v{x}_0, \dots, \v{x}_{k-2})$.
Using Lemma \ref{l:l4} twice, we have
\begin{equation}\label{tmp:12:57_7}
  S \ge g ( \L_{k-1}^* (\v{x}_0, \dots, \v{x}_{k-2}), \v{x}_{k-1} ) - 16 \k \,.
\end{equation}
\\

Let us consider $\ind_{k+1} (\mathbf{\L}_{k+1}) (g)$.
We have
$$
  \ind_{k+1} (\mathbf{\L}_{k+1}) (g) =
$$
$$
               \frac{1}{|\L_0|^2} \sum_{\v{x}_0 \in \L_0} \frac{1}{|\L_1(\v{x}_0)|^2} \sum_{\v{x}_1 \in \L_1 (\v{x}_0)}
               \dots
               \sum_{\v{x}_{k-1} \in \L_{k-1} (\v{x}_0, \dots, \v{x}_{k-2}) } S(\v{x}_0, \dots, \v{x}_{k-1})
$$
By assumption $\ind_{k-1} (\mathbf{\L}_{k-1}, \o{R}_{k-1} ) (g_3) \ge \tau \beta /4$.
In other words
$$
               \frac{1}{|\L_0|^2} \sum_{\v{x}_0 \in \L_0} \frac{1}{|\L_1(\v{x}_0)|^2} \sum_{\v{x}_1 \in \L_1 (\v{x}_0)}
               \dots
$$
\begin{equation}\label{f:tmp_22:17_7}
               \sum_{\v{x}_{k-1} \in \o{R}_{k-1} (\v{x}_0, \dots, \v{x}_{k-2}) }
                    g_3 ( \L_{k-1}^* (\v{x}_0, \dots, \v{x}_{k-2}), \v{x}_{k-1} )
                        \ge
                            \tau \beta /4 \,.
\end{equation}
By $M_{k-1} (\v{x}_0, \dots, \v{x}_{k-2}) $ denote the set of
$\v{x}_{k-1} \in \o{R}_{k-1} (\v{x}_0, \dots, \v{x}_{k-2})$ such that
$\v{x}_{k-1}$ satisfies (\ref{c:g_3}).
Using (\ref{f:tmp_22:17_7}), we obtain
$$
    S_M=
               \frac{1}{|\L_0|^2} \sum_{\v{x}_0 \in \L_0} \frac{1}{|\L_1(\v{x}_0)|^2} \sum_{\v{x}_1 \in \L_1 (\v{x}_0)}
               \dots
$$
\begin{equation}\label{f:tmp_22:24_7}
               \sum_{\v{x}_{k-1} \in M_{k-1} (\v{x}_0, \dots, \v{x}_{k-2}) }
                    g_3 ( \L_{k-1}^* (\v{x}_0, \dots, \v{x}_{k-2}), \v{x}_{k-1} )
                        \ge
                            \tau \beta /8 \,.
\end{equation}
Using (\ref{c:g_3}), (\ref{tmp:12:57`}), (\ref{tmp:12:57_7}) and
(\ref{f:tmp_22:24_7}), we get
$$
  \ind_{k+1} (\mathbf{\L}_{k+1} ) (g)
                \ge
                \frac{1}{|\L_0|^2} \sum_{\v{x}_0 \in \L_0} \frac{1}{|\L_1(\v{x}_0)|^2} \sum_{\v{x}_1 \in \L_1 (\v{x}_0)}
                \dots
$$
$$
                \Big\{
                \sum_{\v{x}_{k-1} \in M_{k-1} (\v{x}_0, \dots, \v{x}_{k-2}) }
                (  g ( \L_{k-1}^* (\v{x}_0, \dots, \v{x}_{k-2}), \v{x}_{k-1} ) ( 1 + 2^{-15} \tau^2 \beta^2 \sigma^3 ) - 8 \k)
                    +
$$
$$
                    +
                \sum_{\v{x}_{k-1} \in \L_{k-1} (\v{x}_0, \dots, \v{x}_{k-2} ) \setminus M_{k-1} (\v{x}_0, \dots, \v{x}_{k-2})}
                (  g ( \L_{k-1}^* (\v{x}_0, \dots, \v{x}_{k-2}), \v{x}_{k-1} ) - 16 \k) \Big\}
             \ge
$$
$$
             \ge
                \ind_{k-1} (\mathbf{\L}_{k-1}) (g)
                +
                  2^{-15} \tau^2 \beta^2 \sigma^3 \Big( \frac{\tau \beta}{8} \Big)
                  S_M
                  - 24 \k
             \ge
$$
$$
                \ge
                    \ind_{k-1} (\mathbf{\L}_{k-1} ) (g) + 2^{-24} \tau^4 \beta^4 \sigma^3  - 24 \k
                \ge
$$
$$
                \ge
                    \ind_{k-1} (\mathbf{\L}_{k-1} ) (g) + 2^{-25} \tau^4 \beta^4 \sigma^3 \,.
$$
In other words, for all $k\ge 1$, we have
\begin{equation}\label{e:inc_ind}
  \ind_{k+1} (\mathbf{\L}_{k+1}) (g) \ge \ind_{k-1} (\mathbf{\L}_{k-1}) (g)
                                                            + 2^{-25} \tau^4 \beta^4 \sigma^3 \,.
\end{equation}

Since for any $k$ we have $\ind_k (\mathbf{\L}_k) (g) \le 1$, it follows that
the total number of steps of the algorithm does not exceed
$K_0=2^{30} \tau^{-4} \beta^{-4} \sigma^{-3}$.

Suppose that the algorithm stops at step $K$, $K\ge 1$,
$K \le 2^{30} \tau^{-4} \beta^{-4} \sigma^{-3}$.
We have
\begin{equation}\label{eq:g_e}
  \ind_{K} (\mathbf{\L}_K, \o{R}_K) (g_3)   <  \frac{\tau \beta}{4} \,.
\end{equation}
Using Lemma \ref{l:keps}, we get
$$
  \ind_{K} (\mathbf{\L}_K) (g_1) \ge (\d + \tau) \beta - 8 \k K  \ge (\d +  \frac{7 \tau}{8} ) \beta \,.
$$
Using (\ref{eq:g_e}), we obtain
\begin{equation}\label{e:t_20:05}
  \ind_{K} ( \mathbf{\L}_K, R_K ) (g_1) \ge (\d +  \frac{3 \tau}{8} ) \beta \,.
\end{equation}
The summation in (\ref{e:t_20:05}) is taken over the sets $\L_K^* (\v{x}_0, \dots, \v{x}_{K-1}) + \v{y}$, where
$\v{y} \in R_K (\v{x}_0, \dots, \v{x}_{K-1}) $.

Let $E_K$ be the family of vectors $\v{y}$ such that
$\v{y} \in E_{K} (\v{x}_0, \dots, \v{x}_{K-1})$, and
$R_K^*$ be the family of vectors  $\v{y}$ such that
$\v{y} \in R_{K} (\v{x}_0, \dots, \v{x}_{K-1})$,
but
$\v{y}$ does not belong to  $E_{K} (\v{x}_0, \dots, \v{x}_{K-1})$.
We have
\begin{equation}\label{tmp:19:47}
  \ind_K (\mathbf{\L}_K, E_K) (g_1)
            < \frac{\tau \beta}{16} \ind_K (\mathbf{\L}_K) (1) \le \frac{\tau \beta}{16}  \,.
\end{equation}
Combining (\ref{e:t_20:05}), (\ref{tmp:19:47}), we get
\begin{equation}\label{tmp:20:25_t}
  \ind_K (\mathbf{\L}_K, R_K^*) (g_1) > (\d + \frac{\tau}{4} ) \beta \,.
\end{equation}
Suppose that for  all $\v{y} \in R_K^* (\v{x}_0, \dots, \v{x}_{K-1})$, we have
$g_2(\L_K^* (\v{x}_0, \dots, \v{x}_{K-1}),\v{y}) < (\d + \tau/16 )$.
Then
$$
  (\d + \frac{\tau}{4} ) \beta  < \ind_K (\mathbf{\L}_K, R_K^*) (g_1)
                        \le (\d + \frac{\tau}{16} ) \ind_{K} (\mathbf{\L}_K, R_K^*) (g_3)
                        \le
$$
\begin{equation}\label{}
                        \le (\d + \frac{\tau}{16} ) \ind_{K} (\mathbf{\L}_K) (g_3) \,.
\end{equation}
Using Lemma \ref{l:keps} once again, we obtain
$$
  (\d + \frac{\tau}{4} ) \beta  < (\d + \frac{\tau}{16} ) \ind_{K} (\mathbf{\L}_K) (g_3)
                        \le (\d + \frac{\tau}{16} ) (\beta + 8 \k K) \le (\d + \frac{\tau}{4} ) \beta
$$
with contradiction.
Whence there exist vectors
$\v{x}_0, \dots, \v{x}_{K-1}$, $\v{y}$ such that
$
  g_2(\L_K^* (\v{x}_0, \dots, \v{x}_{K-1}),\v{y}) \ge (\d + \tau/16 )
$
and
$\v{y} \in R_{K} (\v{x}_0, \dots, \v{x}_{K-1}) \setminus E_{K} (\v{x}_0, \dots, \v{x}_{K-1})$.
Put $\v{t} = \v{y} + \v{s}$ and $\L' = \L_K^* (\v{x}_0, \dots, \v{x}_{K-1})$.
We
obtain the vector $\v{t}$, the sets
$E_1' = (\t{E}_1 - y_1) \cap \L'$, $E_2' = (\t{E}_2 - y_2) \cap \L'$
and the Bohr set $\L'$ which satisfy the conditions $1)$---$3)$.

Let us estimate $D$, $\eps'$ and $N'$.
At the each step of the algorithm the dimension of Bohr sets increases at most $1$.
Since the total number of steps does not exceed $K_0$, it follows that
$D \le d + 2^{30} \tau^{-5} \beta^{-5} \sigma^{-3}$,
$\eps' \ge (2^{-10} \eps)^D \eps_0$, $N' \ge (2^{-10} \eps)^D N$.
Using Lemma \ref{l:Bohr_est} and (\ref{c:Nge}), we obtain that the set  $\L'$ is not empty.
%
%
%
%
This completes the proof.

\refstepcounter{section}

\begin{center}
\textsc{ \arabic{section}. Proof of main result.}
\end{center}

Let us put Theorems \ref{t:Phase1} and \ref{prop:point_3}  together in a single proposition.

\Pred
{ \it
  Let $\L= \L(\theta, \eps_0, N)$ be a Bohr set, $\theta \in {\bf T}^d$,
  and $\v{s} = (s_1,s_2) \in {\bf Z}^2$.
  Let   $E_1$, $E_2$ be sets, $E_i = \beta_i |\L|$, $i=1,2$, $\beta = \beta_1 \beta_2$.
  Suppose
  ${\bf E} = E_1 \m E_2$ is a subset of $(\L+s_1) \m (\L+s_2)$,
  $E_1$, $E_2$ are $(\a_0,2^{-10} \eps^2)$--uniform subsets of $\L + s_1$, $\L + s_2$, respectively,
  $\a_0 = 2^{-2000} \d^{96} \beta_1^{48} \beta_2^{48}$, $\eps = (2^{-100}  \a_0^2 ) / (100 d)$.
  Suppose that $A$ is a subset of ${\bf E}$, $\d_{{\bf E}} (A) = \d$,
  and $A$ has no triples $\{ (k,m), (k+d,m), (k,m+d) \}$ with $d\neq 0$.
  Let
\begin{equation}\label{END!}
    \log N \ge 2^{1000 000} (2^{250 000} \d^{-20 000} \beta^{-200} + d )^3 \log \frac{1}{\d \beta \eps_0 }  \,.
\end{equation}
    Then there is a Bohr set $\t{\L}$
    and a vector $\v{y} = (y_1, y_2) \in {\bf Z}^2$ with the following properties :
    there exist sets $E_1' \subseteq ( E_1  - y_1 \cap \t{\L} ) $, $E_2' \subseteq ( E_2 - y_2 \cap \t{\L} ) $
    such that
\\
\\
$
  \quad 1) ~~
        \mbox{ Let }
        |E_1'| = \beta_1' |\t{\L}|, |E_2'| = \beta_2' |\t{\L}|
                        \mbox{ and } \beta' = \beta_1' \beta_2' \,.
        \mbox{ Then }
         \beta'  \ge 2^{-1500} \d^{100} \beta \,.
$
\\
$
  \quad 2) ~~  E_1', E_2' \mbox{ are } (\a_0', 2^{-10} \eps'^2)\mbox{--uniform, where }
        \a_0' = 2^{-2000} \d^{96} \beta'^{48},
$
\\
$
         \eps' = \frac{ 2^{-100} \a_0'^2}{ 100 D' } \,,
         D\le D' = 2^{250 000} \d^{-20000} \beta^{-200} + d \,.
%
%
$
\\
$
   \quad 3) ~~ \mbox{ For } \t{\L} = \L_{ \t{\theta},\t{\eps},\t{N} } \mbox{ we have }
$
$
   \t{\theta} \in {\bf T}^D ,
$
$
   \t{\eps}
          \ge ( 2^{-100} \eps'^2 )^D  \eps_0
            \mbox{ and }
$
$
    \t{N} \ge (2^{-100} \eps'^2 )^D N \,.
$
\\
$
  \quad 4) ~~
  \d_{E_1' \m E_2' } (A) \ge \d + 2^{-600} \d^{37} \,.
$
}
\label{pred:Ph1_Ph2}
\\

The following lemma is due to B. Green.

\Lemma
{\it
  Let  $N$ be a natural number.
  Suppose  $A$ is a subset of $[-N,N]^2$, $|A| = \d (2N+1)^2$, and
  $A$ has no triples $\{ (k,m), (k+d,m), (k,m+d) \}$ with $d>0$.
  Then there exists a set $A_1 \subseteq A$ such that \\
  $1)~$
  $
   |A_1| \ge \d^2  (2N+1)^2 /4
  $
  and \\
  $2)~$ $ A_1 $ has no triples  $ \{ (k,m), (k+d,m), (k,m+d) \}$ with $d\neq 0$.
}
\label{l:Green's_trick}
\\
\Proof
Since  $|A| = \d (2N+1)^2$, it follows that
\begin{equation}\label{tmp_11:04}
  \sum_{\v{v}} \sum_{\v{s}} A(\v{s}) A(\v{v} - \v{s}) =
  \sum_{\v{s}} \sum_{\v{v}} A(\v{s}) A(\v{v} - \v{s}) = \d^2 (2N+1)^4 \,.
\end{equation}
Clearly, the summation in (\ref{tmp_11:04}) is taken over $\v{v} \in [-2N, -2N+1, \dots, 2N-1, 2N]^2$.
Hence there exists a vector  $\v{v}$ such that
$| A \cap (\v{v} - A) | \ge \d^2 (2N+1)^4 / (4N+1)^2 \ge \d^2  (2N+1)^2 /4$.
Put $A_1 = A \cap (\v{v} - A)$.
We have $A_1 \subseteq A$.
It follows that $A_1$ does not contain any triple  $ \{ (k,m), (k+d,m), (k,m+d) \}$ with $d> 0$.
Since  $A_1 \subseteq \v{v} - A$, it follows that $A_1$ has no triples $ \{ (k,m), (k+d,m), (k,m+d) \}$ with $d< 0$.
This completes the proof.
\\

{\bf Proof of Theorem \ref{main_th}}.

Suppose  $A \subseteq [-N,N]$ and $A$ has no  triples
$\{ (k,m), (k+d,m), (k,m+d) \}$ with $d>0$.
Using Lemma \ref{l:Green's_trick}, we get the set $A'$, $A' \subseteq A$,                           
$|A'| \ge \d^2/4 (2N+1)^2$ such that
$A'$ has no triples
$(k,m), (k+d,m), (k,m+d)$ with $d\neq 0$.
Let $\d' = \d^2 /4$.

The proof of Theorem \ref{main_th} is a sort of an algorithm.

After the $i$th step of the algorithm
an integer vector  $\v{s}_i = (s^{(1)}_i,s^{(2)}_i)$
and sets :
a regular Bohr set $\L_i = \L_{\theta_i, \eps_i, N_i}$,
sets $E_i^{(1)} - s^{(1)}_i \subseteq \L_i $,  $E_i^{(2)} - s^{(2)}_i \subseteq \L_i $,
will be constructed.
Let $|E_i^{(1)}| = \beta^{(1)}_i |\L_i|$,  $|E_i^{(2)}| = \beta^{(2)}_i |\L_i|$,
$\beta_i = \beta^{(1)}_i \beta^{(2)}_i$, ${\bf E_i} = E_i^{(1)} \m E_i^{(2)}$.

The sets $\L_i$, $E_i^{(1)}$, $E_i^{(2)}$ satisfy the following conditions
\\
$ 1)~~~ \beta_{i}  \ge 2^{-1500} \d'^{100} \beta_{i-1} $. \\
$ 2)~~~ E_i^{(1)}, E_i^{(2)}$ are $(\a_0^{(i)}, 2^{-10} (\eps'_{i})^2)$--uniform,
$ \a_0^{(i)} = 2^{-2000} \d'^{96} \beta_i^{48}$, $\eps'_{i} = 2^{-100} (\a_0^{(i)})^2 / (100 d_i )$. \\
$ 3)~~~ \L_i = \L_{\theta_i, \eps_i, N_i},
  \t{\theta} \in {\bf T}^{d_i} , d_i \le 2^{250 000} \d'^{-20000} \beta_{i-1}^{-200} + d_{i-1} ,
$
$
   \eps_{i} \ge (2^{-100} (\eps'_{i})^2 )^{d_i} \eps_{i-1},
$
$
    N_i \ge ( 2^{-100} (\eps'_{i})^2 )^{d_i} N_{i-1}
$.
\\
$ 4)~~~ \d_{\bf E_i} (A') \ge \d_{\bf E_{i-1} } (A') + 2^{-600} \d'^{37} $.
\\

Proposition \ref{pred:Ph1_Ph2} allows us to carry the $(i+1)$th step of the algorithm.
By this Proposition there exists a new vector $\v{s}_{i+1} = (s^{(1)}_{i+1},s^{(2)}_{i+1}) \in {\bf Z}^2$
and sets :
a regular Bohr set
$\L_{i+1} = \L_{\theta_{i+1}, \eps_{i+1}, N_{i+1}}$,
sets
$E_{i+1}^{(1)} - s^{(1)}_i \subseteq \L_{i+1} $,
$E_{i+1}^{(2)} - s^{(2)}_{i+1} \subseteq \L_{i+1} $,
${\bf E_{i+1}} = E_{i+1}^{(1)} \m E_{i+1}^{(2)}$,
which satisfy $1)$ --- $4)$.

Put $\theta_0 = \{ 0 \} $, $\L_0 = \L_{\theta_0, 1, N}$ and $E_1 = E_2 = [-N,N]$, $\beta_0 = 1$.
Clearly,
$E_1$, $ E_2$  are $(2^{-2000} \d'^{96}, 2^{-10 000}  \d'^{400} )$--uniform.
Hence we have constructed zeroth step of the algorithm.

Let us estimate the total number of steps of our procedure.
For an arbitrary $i$ we have $ \d_{\bf E_i} (A') \le 1$.
Using this and condition $4)$, we obtain that the total number of steps cannot be more then
$ 2^{700} \d'^{-36} = K $.

Condition $3)$ implies  $\beta_i \ge (2^{-1500} \d'^{100})^i$.
Hence  $d_i \le (C_1 \d)^{ -C'_1 i }$, where $C_1, C'_1 > 0$ are absolute constants.

To prove  Theorem \ref{main_th}, we                                                                 
need to
verify
condition
(\ref{END!}) at the last step of the algorithm.
Using $3)$, we get
$$
  N_K \ge (C_2 \d )^{C_3 \d^{-C_4 K}} N \,,
$$
where $C_2, C_3, C_4 > 0$ are absolute constants.
Condition (\ref{END!}) can be rewrite as
$$
  N_K \ge (C'_2 \d )^{-C'_3 \d^{-C'_4 K}} \,,
$$
where $C'_2, C'_3, C'_4 > 0$ are absolute constants.
Whence we need to check up the following inequality                                                 
\begin{equation}\label{END!!}
  N \ge (C''_2 \d )^{-C''_3 \d^{-C''_4 K}}  =  \exp ( \d^{ - C' \d'^{-36} }) \,,
\end{equation}
where $C''_2, C''_3, C''_4, C' > 0$ are absolute constants.
By assumption
$$
  \d \gg \frac{1}{ ( \log \log N )^{1/73} } \,.
$$
It follows that
$$
  \d' \gg \frac{1}{ ( \log \log N )^{2/73} }
$$
and we get (\ref{END!!}).
Hence  $A'$ has a triple
$\{ (k,m), (k+d,m), (k,m+d) \}$, where $d\neq 0$.
This contradiction concludes the proof.

\Note Certainly,  the constant $73$ in Theorem  \ref{main_th}  can be slightly decreased.
Nevertheless, it is the author's opinion that this constant cannot be lowered  as to such $1$
without a new idea.



\refstepcounter{section}


\begin{center}
\textsc{ \arabic{section}. On quantitative  recurrence.}
\end{center}

In this section we apply  Theorem \ref{main_th} to the theory of dynamical systems.

Let $X$ be a metric space with metric $d(\cdot ,\cdot )$ and a
Borel sigma--algebra of measurable sets $\Phi$. Let $T$ be a
measure preserving transformation of a measure space
$(X,\Phi,\mu)$, and let us assume that measure of $X$ is equal to
$1$. The well--known Poincare  theorem (see \cite{Poin}) asserts
that for almost every point $x \in X $:
$$ \forall \eps >  0 ~\forall K  >  0 ~\exists
t  >  K : d(T^tx,x) < \eps .$$
Consider a measure $H_h (\cdot)$ on
$X$, defined as
$$
  H_h (E)= \lim_{\delta \rightarrow 0} H_{h}^{\delta }(E),
$$
where $h(t)$ is a positive ($ h(0)=0 $) continuous increasing
function and $H_{h}^{\delta }(E)= \inf_\tau \{ \sum h(\delta_j)
\}$, when $\tau$ runs through all countable coverings $E$ by open
sets
$ \{ B_j \}$ , $diam(B_i) = \delta_j < \delta$.\\
If $h(t) = t^\alpha$, then we get the ordinary Hausdorff measure
$H_{\alpha } (\cdot ) $.

We shall say that a measure $\mu$ is congruent to a measure $H_h$, if
any $\mu$--measurable set is $H_h$--measurable.

The following theorems
\arabic{section}.1 and \arabic{section}.4
were proven in
\cite{Sh} (see also \cite{Bo,Mo}).

\Th
{\it Let {\it X } be a metric space with
$ H_h (X) = C < \infty $,
and {\it T } be a measure preserving transformation of {\it X}.
Assume that $\mu$ is congruent to $H_h$.
\\Consider the following function:
$
 C(x) = \liminf_{n \to \infty}  \{ n \cdot h(d(T^{n}x,x)) \}
$.
\\
Then the function $C(x)$ is $\mu$--integrable
and for any $\mu$--measurable set $A$, we have
$$
  \int_A C(x) d\mu \le H_h (A).
$$
If $H_h(A)=0$, then $\int_A C(x) d\mu = 0$
with no demand on measures $\mu$ and $H_h$ to be congruent.
}
\label{x_1}

Now we introduce the following concept (see \cite{Kolm}).

\Def
 Let $G$ be a totally bounded subset of $X$.
 By $N_{\eps}(G,X)$ denote the minimal cardinality of $\eps$--net of G.
Put $N_\eps (X) = N_\eps (X,X)$.
\\
If  $X$ is totally bounded, then for any $\delta$, we have
$N_\delta (X)< \infty$ and $\sum h(\delta_j) \le N_{\delta}(X) h(\delta )$.
Let $h$ be the function from the definition of
$H_h$. If $ N_{\delta }(X) \le C/h(\delta )$, then $H_{h}(X) \le C$.

\Def Let $N$ be a natural number.
By $C_{N} (x)$ denote the function
$C_{N}(x) = \min \{~ d(T^{n}x,x) ~|~ 1 \le n \le N ~\}$.
\label{const_rec}

\Th {\it Let $X$ be a totally bounded metric space with metric
$d(\cdot ,\cdot )$, and let $N(x) = N_x (X)$.
Suppose that ${\it T}$ is a measure--preserving
transformation of {\it X}, and $diam (X) = 1$.
\\
Let $A \subseteq X$ be an arbitrary $\mu$--measurable set, and
let $g(x)$ be the real nondecreasing
function bounded on $[0,1]$ such that
for any $t \in (0,1]$ there exists Stieltjes
integral $\int_t^{1} N_{A}(x) dg(x)$,
where $N_{A}(x) = \min (\mu (A),N_x(A,X)/N)$.
Then
$$
  \int_{A} g(C_{N}(x)) d\mu \le
  \inf_t \{ g(t)\mu(A) +
  \int_t^1 N_{A}(x) dg(x) ~\}.
$$
} \label{x_2}

The following lemma  is due to Poincare (see \cite{Poin,Bo}).
\\
\Lemma {\it Let $Y$ be a $\mu$--measurable set, and  $ t \ge 1 $.
Set
$$
  Y(t) := \{ x \in Y ~|~ T^{i}x \notin Y
  \mbox{ for all natural } i, 1 \le i \le t \} .
$$
Then $ \mu(Y(t)) \le 1/t .$
}\label{l_add}

  This lemma is the main tool in proving  Theorems \ref{x_1}, \ref{x_2}.
  \\
  %
  %
  %
  %

  Let us now consider the case of two commutative operators.
  Let
$S$ and $R$ be two {\it commutative} measure--preserving transformation of $X$.
The next result is the main one of this section.

\Th
{\it Let {\it X } be a metric space with
$ H_h (X) = C < \infty $,
and let  $S,R$
be two  commutative measure--preserving transformation of $X$.
Assume that $\mu$ is congruent to $H_h$.
\\
Let us consider the function
$$
 C_{S,R}(x) =
 \liminf_{n \to \infty}  \{L^{-1}(n) \cdot max \{ h(d(S^{n}x,x)), h(d(R^{n}x,x)) \} \},
$$
where $L^{-1}(n) = 1/ L(n)$.
\\
Then the function $C_{S,R}(x)$ is $\mu$--integrable
and for any $\mu$--measurable set $A$, we have
$$
  \int_A C_{S,R}(x) d\mu \le H_h (A).
$$
If $H_h(A)=0$, then $\int_A C_{S,R}(x) d\mu = 0$
with no demand on measures $\mu$ and $H_h$ to be congruent. } \label{x_3}

The next definition is analog of Definition \ref{const_rec}.

\Def Let $N$ be a natural number.
By $C^{S,R}_{N} (x)$ denote the function
$C^{S,R}_{N}(x) = \min \{~ max \{ d(S^{n}x,x), d(R^{n}x,x) \} ~|~ 1 \le n \le N ~\}$.
The function $C^{S,R}_{N} (x)$ will be called
{\it $N$--constant of simultaneously recurrence} for point $x$.

\Th {\it Let $X$ be a totally bounded metric space
with metric
$d(\cdot ,\cdot )$, and $N(x) = N_x (X)$.
Suppose that
${\it S,R}$ are two measure--preserving transformation of {\it X},
and $diam (X)= 1$.
\\
Let $A \subseteq X$ be an arbitrary  $\mu$--measurable set, and
let $g(x)$
be the real nondecreasing
function bounded on $[0,1]$ such that
for any $t \in (0,1]$ there exists Stieltjes
integral $\int_t^{1} N_{A}(x) dg(x)$,
where
 $N_{A}(x) = \min (\mu (A),N_x(A,X) L(N))$.
Then
$$
  \int_{A} g(C^{S,R}_N(x)) d\mu \le
  \inf_t \{ g(t)\mu(A) +
  \int_t^1 N_{A}(x) dg(x) ~\}.
$$
}
\label{x_4}

%
%
%
%


The next Lemma is the main of this section.
Using this lemma we obtain Theorems
\ref{x_3} and \ref{x_4} by the same argument as Lemma \ref{l_add} implies Theorems
\ref{x_1} and \ref{x_2}
(for details see \cite{Sh}).

\Lemma {\it Let $Y$ be a $\mu$--measurable set, $ t \ge 1 $.
Set
$$
  Y(t) := \{ x \in Y ~|~ \mbox{ Either } S^{i}x \notin Y \mbox{ or } R^{i}x \notin Y
  \mbox{ for all natural } i, 1 \le i \le t \} .
$$\label{l_key}
Then
$\mu(Y(t)) \le L(t)$.
}
\label{ll}
\\
\Proof
See \cite{Shkr_mult}.

 Now we apply Theorem \ref{x_3} to the case of compact metric space.

The following lemma can be found in \cite{Bg}.

\Lemma {\it Let $X$ be a compact metric space, and let
$T_1, \dots T_l$ be
continuous commutative transformations of $X$.
Then there exists a finite measure $\mu$
such that transformations $T_1, \dots T_l$ preserve $\mu$.}

\Cor {\it Let $X$  be a compact metric space with metric $d(\cdot,\cdot)$,
and
$H_h (X) < \infty$.
Let $~S, R$ be two continuous commutative transformations of $X$.
Then there exists  $x\in X$ such that
\\
$
 \liminf_{n \to \infty}  \{L^{-1}(n) \cdot max \{ h(d(S^{n}x,x)), h(d(R^{n}x,x)) \} \}
 \le H_h (X).
$
}

\end{document}